\documentclass[english, a4paper, article, twoside, 10pt, leqno]{memoir}

\usepackage[draft]{fixme}

\usepackage[T1]{fontenc} 	%European fonts
\usepackage[noDcommand,slantedGreeks]{kpfonts}

\usepackage{babel}
\usepackage{microtype}

\usepackage[shortlabels]{enumitem}

\raggedbottomsectiontrue%less harse than \raggedbottom

\makeatletter

\AddToHook{bfseries}{\boldmath}
\usepackage[autostyle,english=american]{csquotes}

\usepackage[backend=biber,citestyle=authoryear-icomp,bibstyle=authoryear,
	uniquename=init,giveninits=true,autolang=hyphen]{biblatex}
\nobibintoc %% Removes bibliography from the pable of contents.

\DeclareDocumentCommand\textdef{m}{% for inline text definitions
	\textbf{#1}%
}

\usepackage[thmmarks,amsmath]{ntheorem}

\newif\ifprenumber
\prenumbertrue

\makeatletter

\makeatletter
\newtheoremstyle{notitle}% for exercises
{\item[\hskip\labelsep \theorem@headerfont ##2\theorem@separator]}%
{\item[\hskip\labelsep \theorem@headerfont ##2\theorem@separator\ ##3\theorem@separator]}

\newtheoremstyle{notitlebreak}% for exercises with a break
{\item[\rlap{\vbox{\hbox{\hskip\labelsep \theorem@headerfont ##2\theorem@separator}\hbox{\strut}}}]}%
{\item[\rlap{\vbox{\hbox{\hskip\labelsep \theorem@headerfont ##2\theorem@separator\ ##3\theorem@separator}\hbox{\strut}}}]}

\newtheoremstyle{prenumber}%
  {\item[\hskip\labelsep \theorem@headerfont ##2\theorem@separator\ ##1\theorem@separator]}%
  {\item[\hskip\labelsep \theorem@headerfont ##2\theorem@separator\ ##3\theorem@separator]}

\newtheoremstyle{prenumberbreak}%
  {\item[\rlap{\vbox{\hbox{\hskip\labelsep \theorem@headerfont 
          ##2\theorem@separator\ ##1\theorem@separator}\hbox{\strut}}}]}%
  {\item[\rlap{\vbox{\hbox{\hskip\labelsep \theorem@headerfont 
          ##2\theorem@separator\ ##3\theorem@separator}\hbox{\strut}}}]}

\newtheoremstyle{customplain}%
  {\item[\hskip\labelsep \theorem@headerfont ##1\ ##2\theorem@separator]}%
  {\item[\hskip\labelsep \theorem@headerfont ##3\ ##2\theorem@separator]}
  
\newtheoremstyle{custombreak}%
  {\item[\rlap{\vbox{\hbox{\hskip\labelsep \theorem@headerfont
          ##1\ ##2\theorem@separator}\hbox{\strut}}}]}%
  {\item[\rlap{\vbox{\hbox{\hskip\labelsep \theorem@headerfont 
          ##3\ ##2\theorem@separator}\hbox{\strut}}}]}

\newtheoremstyle{postnumpara}%
  {\item[\hskip\labelsep \theorem@headerfont ##2\theorem@separator]}%
  {\item[\hskip\labelsep \theorem@headerfont ##3\ ##2\theorem@separator]}
  
\newtheoremstyle{postnumparabreak}%
  {\item[\rlap{\vbox{\hbox{\hskip\labelsep \theorem@headerfont
          ##2\theorem@separator}\hbox{\strut}}}]}%
  {\item[\rlap{\vbox{\hbox{\hskip\labelsep \theorem@headerfont 
          ##3\ ##2\theorem@separator}\hbox{\strut}}}]}

\newtheoremstyle{nonumberplainflex}% for proof environments, see below
	{\item[\theorem@headerfont\hskip\labelsep ##1\theorem@separator]}%
	{\item[\theorem@headerfont\hskip \labelsep ##3\theorem@separator]}

\newtheoremstyle{nonumberbreakflex}%
	{\item[\rlap{\vbox{\hbox{\hskip\labelsep \theorem@headerfont
		##1\theorem@separator}\hbox{\strut}}}]}%
	{\item[\rlap{\vbox{\hbox{\hskip\labelsep \theorem@headerfont 
		##3\theorem@separator}\hbox{\strut}}}]}

\theoremseparator{.}

\theorembodyfont{\normalfont}

\ifprenumber
	\theoremstyle{notitle}
\else
	\theoremstyle{postnumpara}
\fi

\ifprenumber
	\theoremstyle{prenumber}
\else
	\theoremstyle{customplain}
\fi

\theorembodyfont{\itshape}

\newtheorem{theorem}[equation]{Theorem}
\newtheorem{corollary}[equation]{Corollary}
\newtheorem{proposition}[equation]{Proposition}
\newtheorem{lemma}[equation]{Lemma}
\newtheorem{conjecture}[equation]{Conjecture}

\theorembodyfont{\normalfont}
\newtheorem{definition}[equation]{Definition}
\newtheorem{notation}[equation]{Notation}
\theoremsymbol{\ensuremath{\triangle}}
\newtheorem{remark}[equation]{Remark}

\theoremsymbol{\ensuremath{\scriptstyle\bigcirc}}
\newtheorem{example}[equation]{Example}

\theoremstyle{nonumberplainflex}
\theoremsymbol{\ensuremath{\square}}
\theoremheaderfont{\itshape\bfseries}
\newtheorem{proof}{Proof}

\theoremstyle{nonumberplainflex}
\theoremheaderfont{\normalfont\bfseries}
\theoremsymbol{}

\theorembodyfont{\itshape}

\theorembodyfont{\normalfont}

\theoremsymbol{\ensuremath{\triangle}}

\theoremsymbol{\ensuremath{\scriptstyle\bigcirc}}

\theoremsymbol{}
\theoremheaderfont{\normalfont\bfseries}

\ifprenumber
	\theoremstyle{prenumberbreak}
\else
	\theoremstyle{custombreak}
\fi

\theorembodyfont{\itshape}

\theorembodyfont{\normalfont}
\newtheorem{definitionbreak}[equation]{Definition}
\theoremsymbol{\ensuremath{\triangle}}

\theoremsymbol{\ensuremath{\scriptstyle\bigcirc}}

\theoremsymbol{}

\theoremstyle{nonumberbreakflex}
\theoremsymbol{\ensuremath{\square}}
\theoremheaderfont{\itshape\bfseries}

\theoremstyle{nonumberbreakflex}
\theoremheaderfont{\normalfont\bfseries}
\theoremsymbol{}
\theorembodyfont{\itshape}

\theorembodyfont{\normalfont}

\theoremsymbol{\ensuremath{\triangle}}

\theoremsymbol{\ensuremath{\scriptstyle\bigcirc}}

\usepackage{tikz}
\usetikzlibrary{babel}
\usetikzlibrary{cd}

\tikzcdset{
  arrow style=tikz,
}

\tikzset{
  >/.tip={Stealth[length=2.9pt, width=4.4pt, inset=1.8pt]},
  tikzcd left hook/.tip={Hooks[
	  left,
	  length=2pt,
	  width=5.5pt,
	]},
  iso/.style={
    every to/.append style={
      edge node={
        node [sloped, allow upside down]{
          \raisebox{0.1em}[0pt][0pt]{\ensuremath{\sim}}
        }
      }
    }
  },
  iso'/.style={
    every to/.append style={
      edge node={
        node [sloped, allow upside down]{
          \raisebox{-0.6em}[0pt][0pt]{\ensuremath{\sim}}
        }
      }
    }
  },
  symbol/.style={
      draw=none,
      every to/.append style={
        edge node={node [sloped, allow upside down, auto=false]{$#1$}}}
  },
}

\usepackage[hidelinks]{hyperref}

\newif\ifzcrefclever
\zcrefclevertrue

\ifzcrefclever
	\usepackage{zref-clever}
	
	\zcRefTypeSetup{theorem}{
		Name-sg = Theorem,
		name-sg = Theorem,
		Name-pl = Theorems,
		name-pl = Theorems,
	}
	\zcRefTypeSetup{definition}{
		Name-sg = Definition,
		name-sg = Definition,
		Name-pl = Definitions,
		name-pl = Definitions,
	}
	\zcRefTypeSetup{example}{
		Name-sg = Example,
		name-sg = Example,
		Name-pl = Examples,
		name-pl = Examples,
	}
	\zcRefTypeSetup{corollary}{
		Name-sg = Corollary,
		name-sg = Corollary,
		Name-pl = Corollaries,
		name-pl = Corollaries,
	}
	\zcRefTypeSetup{proposition}{
		Name-sg = Proposition,
		name-sg = Proposition,
		Name-pl = Propositions,
		name-pl = Propositions,
	}
	\zcRefTypeSetup{lemma}{
		Name-sg = Lemma,
		name-sg = Lemma,
		Name-pl = Lemmas,
		name-pl = Lemmas,
	}
	\zcRefTypeSetup{remark}{
		Name-sg = Remark,
		name-sg = Remark,
		Name-pl = Remarks,
		name-pl = Remarks,
	}
	\zcRefTypeSetup{conjecture}{
		Name-sg = Conjecture,
		name-sg = Conjecture,
		Name-pl = Conjectures,
		name-pl = Conjectures,
	}
	
	\newcommand\@add@theorem@hooks[1]{%
		\AddToHook{env/#1/begin}{%
			\zcsetup{countertype={equation=#1}}}
		\AddToHook{env/#1break/begin}{%
			\zcsetup{countertype={equation=#1}}}
	}

	\@add@theorem@hooks{theorem}
	\@add@theorem@hooks{definition}
	\@add@theorem@hooks{example}
	\@add@theorem@hooks{corollary}
	\@add@theorem@hooks{proposition}
	\@add@theorem@hooks{lemma}
	\@add@theorem@hooks{remark}
	\@add@theorem@hooks{conjecture}
\else
	\usepackage{cleveref}
	
	\crefname{theorem}{Theorem}{Theorems}
	\crefname{definition}{Definition}{Definitions}
	\crefname{example}{Example}{Examples}
	\crefname{corollary}{Corollary}{Corollaries}
	\crefname{proposition}{Proposition}{Propositions}
	\crefname{lemma}{Lemma}{Lemmata}
	\crefname{remark}{Remark}{Remarks}
	\crefname{claim}{Claim}{Claims}
	\crefname{exercise}{Exercise}{Exercises}
	\crefname{textexercise}{Exercise}{Ecercises}
	\crefname{conjecture}{Conjecture}{Conjectures}
	\crefname{fact}{Fact}{Facts}

	\crefalias{theorembreak}{theorem}
	\crefalias{definitionbreak}{definition}
	\crefalias{examplebreak}{example}
	\crefalias{corollarybreak}{corollary}
	\crefalias{propositionbreak}{proposition}
	\crefalias{lemmabreak}{lemma}
	\crefalias{claimbreak}{claim}
	\crefalias{remarkbreak}{remark}	
	\crefalias{conjecturebreak}{conjecture}
	
	\NewDocumentCommand\zcref{om}{%
		\cref{#2}%
	}
\fi

\newlist{definitionlist}{enumerate}{2}
\setlist[definitionlist]{
	label=\textup{(\roman*)},
	ref=\theequation\textup{(\roman*)},
	resume,
}

\def\varlim@@auxiliary#1#2#3{%
  \vtop{\m@th\ialign{##\cr
    \hfil$#1\operator@font #3$\hfil\cr
    \noalign{\nointerlineskip\kern0\ex@}
    \expandafter#2\ifx#1\scriptscriptstyle\scriptscriptstyle\else\scriptstyle\fi\cr
    \noalign{\nointerlineskip\kern-\ex@}\cr}}%
}
\def\varlim@@#1#2{%
	\varlim@@auxiliary#1#2%
}

\newcommand\dirlim@format[1]{\mathop{\mathpalette\varlim@@{{\rightarrowfill@}{#1}}}\nmlimits@}
\newcommand\invlim@format[1]{\mathop{\mathpalette\varlim@@{{\leftarrowfill@}{#1}}}\nmlimits@}

\NewDocumentCommand\dirlimcommand{me{_^}}{%
	\mathchoice
		{ \dirlim@format { #1 } \IfValueT{#2} { _{#2} } \IfValueT {#3} { {} ^ { #3 } } }
	    { \dirlim@format { #1 } \IfValueT{#2} { _{#2} } \IfValueT {#3} { ^ { #3 } } }
	    { \dirlim@format { #1 } \IfValueT{#2} { _{#2} } \IfValueT {#3} { ^ { #3 } } }
	    { \dirlim@format { #1 } \IfValueT{#2} { _{#2} } \IfValueT {#3} { ^ { #3 } } }
}

\NewDocumentCommand\invlimcommand{me{_^}}{%
	\mathchoice
		{ \invlim@format { #1 } \IfValueT{#2} { _{#2} } \IfValueT {#3} { {} ^ { #3 } } }
	    { \invlim@format { #1 } \IfValueT{#2} { _{#2} } \IfValueT {#3} { ^ { #3 } } }
	    { \invlim@format { #1 } \IfValueT{#2} { _{#2} } \IfValueT {#3} { ^ { #3 } } }
	    { \invlim@format { #1 } \IfValueT{#2} { _{#2} } \IfValueT {#3} { ^ { #3 } } }
}

\newcommand\diff[1]{\mathop{}\!#1}
\def\differentialdspace{\mskip-.7\thinmuskip}
\def\differentialdeltaspace{\mskip-.3\thinmuskip}

\newcommand\wherecommand[1]{%
 \nonscript\:%
 #1\vert
 \allowbreak
 \nonscript\:%
 \mathopen{}%
 }

\newcommand\parenthesesright[1]{%
	\mathrlap{%
		\qquad\qquad\qquad\qquad
		\qquad\qquad
		#1%
	}%
}

\let\save@mathaccent\mathaccent
\newcommand*\if@single[3]{%
  \setbox0\hbox{${\mathaccent"0362{#1}}^H$}%
  \setbox2\hbox{${\mathaccent"0362{\kern0pt#1}}^H$}%
  \ifdim\ht0=\ht2 #3\else #2\fi
  }
\newcommand*\rel@kern[1]{\kern#1\dimexpr\macc@kerna}
\newcommand\widebar{}% initialize
\DeclareRobustCommand*\widebar[1]{\@ifnextchar^{{\wide@bar{#1}{0}}}{\wide@bar{#1}{1}}}
\newcommand*\wide@bar[2]{\if@single{#1}{\wide@bar@{#1}{#2}{1}}{\wide@bar@{#1}{#2}{2}}}
\newcommand*\wide@bar@[3]{%
  \begingroup
  \def\mathaccent##1##2{%
    \let\mathaccent\save@mathaccent
    \if#32 \let\macc@nucleus\first@char \fi
    \setbox\z@\hbox{$\macc@style{\macc@nucleus}_{}$}%
    \setbox\tw@\hbox{$\macc@style{\macc@nucleus}{}_{}$}%
    \dimen@\wd\tw@
    \advance\dimen@-\wd\z@
    \divide\dimen@ 3
    \@tempdima\wd\tw@
    \advance\@tempdima-\scriptspace
    \divide\@tempdima 10
    \advance\dimen@-\@tempdima
    \ifdim\dimen@>\z@ \dimen@0pt\fi
    \rel@kern{0.6}\kern-\dimen@
    \if#31
      \overline{\rel@kern{-0.6}\kern\dimen@\macc@nucleus\rel@kern{0.4}\kern\dimen@}%
      \advance\dimen@0.4\dimexpr\macc@kerna
      \let\final@kern#2%
      \ifdim\dimen@<\z@ \let\final@kern1\fi
      \if\final@kern1 \kern-\dimen@\fi
    \else
      \overline{\rel@kern{-0.6}\kern\dimen@#1}%
    \fi
  }%
  \macc@depth\@ne
  \let\math@bgroup\@empty \let\math@egroup\macc@set@skewchar
  \mathsurround\z@ \frozen@everymath{\mathgroup\macc@group\relax}%
  \macc@set@skewchar\relax
  \let\mathaccentV\macc@nested@a
  \if#31
    \macc@nested@a\relax111{#1}%
  \else
    \def\gobble@till@marker##1\endmarker{}%
    \futurelet\first@char\gobble@till@marker#1\endmarker
    \ifcat\noexpand\first@char A\else
      \def\first@char{}%
    \fi
    \macc@nested@a\relax111{\first@char}%
  \fi
  \endgroup
}

\numberwithin{equation}{chapter}

\newcommand\ABBRVMNrule{Murnaghan--Nakayama rule}
\newcommand\ABBRVlrcoef{Littlewood--Richardson coefficient}

\usepackage{ytableau}

\ytableausetup{centertableaux}

\title{The ring of \( \omega \)-invariant symmetric functions in characteristic~\( 2 \)}

\date{\today}
\author{Sebastian Ørsted}

\hypersetup{
	pdfauthor={\@author},
	pdftitle={The ring of ω-invariant symmetric functions in characteristic 2},
	pdfcreator={LaTeX},
}

\usepackage{mathtools,semtex}

\begin{document}

\maketitle

\begin{abstract}
	\noindent
	We provide a simple presentation by generators and relations
	of the ring of \( \omega \)-invariant symmetric functions
	over the field~\( \mathbb{F}_{2} \).
	Here, \( \omega \)~denotes the standard involution
	on the ring of symmetric functions,
	interchanging the elementary symmetric functions with the complete homogeneous symmetric functions.
	Along the way, we prove several important properties
	of this involution in the specific setting of characteristic~\( 2 \).
	
	\medskip
	
	invariants, Schur functions, Young diagrams, Grassmannians
	
	\medskip
	
\end{abstract}

\noindent
We consider the ring~\( R \) of symmetric functions
over the field~\( \mathbb{F}_{2} \).
This can be defined
\parencite[see e.g.][section~6.2]{ful}
as the direct limit~\(
	R
	=
	\dirlimcommand{lim} R_{(n)} \)
of the rings~\(
	R_{(n)}
	\subset
	\mathbb{F}_{2} \lbrack x_{1}, \dotsc, x_{n} \rbrack
\)
of symmetric polynomials in \( n \)~variables.
This~\( R_{(n)} \) is itself a polynomial ring
in \( n \)~variables, given either by the \emph{elementary symmetric polynomials}
\begin{align}
	e_{k}
	&
	\parenthesesright{( 1 \le k \le n )}
	=
	\sum_{1\le i_{1}<\dotsb <i_{k} \le n} x_{i_{1}} \dotsm x_{i_{k}} \nonumber
\\
\intertext{or the \emph{complete homogeneous symmetric polynomials}}
	h_{k}
	&
	\parenthesesright{( 1 \le k \le n )\text{.}}
	=
	\sum_{1\le i_{1} \le \dotsb \le i_{k} \le n} x_{i_{1}} \dotsm x_{i_{k}} \nonumber
\\
\intertext{The two are related by the formula}
	0
	&
	\parenthesesright{( 1 \le k \le n )\text{,}}
	=
	\sum_{i =0}^{k} e_{i} h_{k -i} 	\label{eq:e_h_relation}
\end{align}
where we use the convention that~\( e_{0} = h_{0} = 1 \).

By extension,
\( R \)~also becomes a polynomial algebra in infinitely many variables,
given by either \( e_{k} \) or~\( h_{k} \).
One of the traditional ways of describing the ring~\( R \)
is as linear combinations of the
\emph{Schur functions}~\( s_{\lambda} \)
with respect to partitions~\( \lambda \).
This association is described in detail in \textcite{ful}.
In particular, taking~\( \lambda = (n) \), we
have~\(
	s_{\lambda} = h_{n}
\),
while for~\( \lambda = (1^{n}) \), the diagram consisting of~\( n \)~copies of~\( 1 \),
we have~\( s_{\lambda} = e_{n} \).
The product structure is given by~\smash{\(
	s_{\lambda} s_{\mu} 	=
	\sum_{\nu} c_{\lambda, \mu}^{\nu} \, s_{\nu} \)},
where \smash{\( c_{\lambda, \mu}^{\nu} \)}~are
the \emph{\ABBRVlrcoef{}s},
which have a non-trivial definition.
They are symmetric in the sense that~\smash{\(
	c_{\lambda, \mu}^{\nu}
	=
	c_{\mu, \lambda}^{\nu}
\)},
and we have
\smash{\(
	c_{\lambda, \mu}^{\nu}
	=
	0
\)}
unless \( \nu \)~contains both
\( \lambda \) and~\( \mu \),
and the sizes of the partitions are related by~\(
	\lvert \lambda \rvert
	+
	\lvert \mu \rvert
	=
	\lvert \nu \rvert
\).

The assignment~\( e_{k} \mapsto h_{k} \)
defines an involution on~\( R \) commonly denoted by~\( \omega \colon R \to R \).
For brevity, we shall usually write it as~\( x \mapsto \widebar{x} \)
instead and pretend that it is complex conjugation.
The \( \omega \)-involution
has a neat description in terms of Schur functions as~\( \widebar{s}_{\lambda} = s_{\lambda^{\!\vee}} \), where \( \lambda^{\!\vee} \)~denotes the conjugate partition of~\( \lambda \).
If a partition is visualized by its Young diagram, the conjugate partition is given by the mirrored diagram.

Because of this description, the subring~\( S = R^{\omega} \)
of \( \omega \)-invariants
is equal to the \( \mathbb{F}_{2} \)-span
of two different classes of vectors:
on the one hand, \( s_{\lambda} \)~for self-conjugate partitions~\( \lambda \),
i.e.~partitions with~\( \lambda = \lambda^{\!\vee} \);
on the other hand,
\( s_{\lambda} + s_{\lambda^{\!\vee}} \)~for~\( \lambda \)
with~\( \lambda \neq \lambda^{\!\vee} \).
That, in principle, fully determines the ring~\( S \),
but this description can be difficult to work with due to the complexity of calculating the \ABBRVlrcoef{}s.
The main goal of this article is to provide a presentation by generators and relations for~\( S \),
which will happen in \zcref{chap:generators_and_relations_for_S}.
To my knowledge, this presentation, as well as the other main results in this paper,
have not appeared in the literature before.

Since \( R \)~is a polynomial ring in both
\( e_{k} \) and~\( h_{k} \),
we shall mostly discard the ground variables~\( x_{k} \)
and simply work with a polynomial ring
\[
	R
	=
	\mathbb{F}_{2} \lbrack w_{1}, w_{2}, w_{3}, \dotsc \rbrack
\]
without specifying whether we realize~\( w_{i} \)
as \( e_{i} \) or~\( h_{i} \).
In this case, the equation~\eqref{eq:e_h_relation}
can be taken as a recursive definition of the involution~\( \omega \).

The intended applications are topological in nature,
arising from the classical fact that
the \( \mathbb{F}_{2} \)-cohomology of the real Grassmannian manifold~\( \mathup{Gr} (n, m) \)
of all \( n \)-planes inside~\( \mathbb{R}^{n +m} \)
is given by
\[
	H^{\SemantexBullet} (\mathup{Gr} (n, m);\mathbb{F}_{2})
	=
	\mathbb{F}_{2} \lbrack w_{1}, \dotsc, w_{n} \rbrack /(\widebar{w}_{m +1}, \dotsc, \widebar{w}_{m +n})
\text{,}
\]
where~\(
	w_{i} = w_{i} (\gamma)
\)
denotes the Stiefel--Whitney class of the canonical \( n \)-plane
bundle~\( \gamma \) over~\(
	\mathup{Gr} (n, m)
\)
\parencite[see e.g.][Proposition~11.1]{borel},
and where \( \widebar{w}_{k} \)~is interpreted using the same formulae as in the ring~\( R \), with the understanding that~\( w_{j} = 0 \) for~\( j > n \).
In the case~\( m = n \), the involution~\( \omega \) on~\(
	H^{\SemantexBullet} (\mathup{Gr} (n, n);\mathbb{F}_{2})
\)
corresponds to the geometric operation of
sending an \( n \)-plane~\smash{\(
	V
	\subset
	\mathbb{R}^{2n }
\)}
to its orthogonal complement~\(
	V^{\perp}
	\subset
	\mathbb{R}^{2n }
\).
The present paper is part of a larger project aimed at calculating the \( \mathbb{F}_{2} \)-cohomology of the projective Grassmannian
manifold, obtained from the ordinary Grassmannian
by identifying~\( V \)
with~\( V^{\perp} \).
Therefore, in \zcref{chap:finite_grassmannian}, we use our presentation of~\( R^{\omega} \) to derive a presentation
of the ring of \( \omega \)-invariants
in~\(
	H^{\SemantexBullet} (\mathup{Gr} (n, n);\mathbb{F}_{2})
\).

\section*{Acknowledgments}

Special thanks to Marcel Bökstedt for suggesting the problem and for many fruitful discussions, helpful comments, and valuable suggestions, as well as for being an excellent NYT~\emph{Connections} partner.
Thanks also to Robert Bruner
for insightful discussions and for providing machine calculations that partly inspired the generators-and-relations presentation in \zcref{chap:generators_and_relations_for_S}.

\newpage

\chapter{Thick differentials, transversality, and normality}\label{chap:general_ring_vocabulary}

In this chapter, and this chapter only, \( R \)~denotes an arbitrary commutative \( \mathbb{F}_{2} \)-algebra with an involution~\( \omega \colon R \to R \) which we write as~\( x \mapsto \widebar{x} \) for brevity.
We consider the subalgebra~\( S = R^{\omega} \)
of \( \omega \)-invariants as well as the ideal
\[
	I
	=
	\lbrace \, x +\widebar{x} \wherecommand{} x \in R \, \rbrace
	\subset
	S
	\text{.}
\]
We denote by~\( \diff{d} \colon R \to R \)
the \( S \)-linear map~\( \diff{d\differentialdspace} x = x + \widebar{x} \)
and note that \( S \)~is the kernel of~\( \diff{d} \) and~\( I \) the image.
As the notation suggests, we shall think of~\( \diff{d} \)
as a kind of differential. Indeed, we have~\( \diff{d}^{2} = 0\),
and \( \diff{d} \)~satisfies the \enquote{Thick Leibniz Rule}
\begin{equation}\label{eq:thick_leibniz}
	\diff{d} (x y )
	=
	x \diff{d\differentialdspace} y 	+
	y \diff{d\differentialdspace} x 	+
	\diff{d\differentialdspace} x \diff{d\differentialdspace} y 	\text{.}
\end{equation}
In the case of three variables, the pattern of the rule
becomes more apparent:
\begin{align*}
	\diff{d} (x y z )
		&= x y \diff{d\differentialdspace} z 		+ x z \diff{d\differentialdspace} y 		+ y z \diff{d\differentialdspace} x \\
		&\quad{}+ x \diff{d\differentialdspace} y \diff{d\differentialdspace} z 		+ y \diff{d\differentialdspace} x \diff{d\differentialdspace} z 		+ z \diff{d\differentialdspace} x \diff{d\differentialdspace} y \\
		&\quad{} + \diff{d\differentialdspace} x \diff{d\differentialdspace} y \diff{d\differentialdspace} z 		\text{.}
\end{align*}
For a general number of elements, the rule becomes
\begin{equation}\label{eq:thick_leibniz_general}
    \diff{d} (x_{1} \dotsm x_{n} )
    =
    \sum_{\varnothing \neq T \subset \lbrace 1, \dotsc, n \rbrace} x_{T^{\complement}} \, (\diff{d\differentialdspace} x )_{T}     \text{.}
\end{equation}
Here, \( (\diff{d\differentialdspace} x )_{T} \)~means the product of~\( \diff{d\differentialdspace} x_{i} \) for~\( i \in T \) while~\(  x_{T^{\complement}} \) means the product of~\( x_{j} \) for~\smash{\( j \in T^{\complement} \)}.

\begin{definition}
	The involution~\( \omega \)
	is called \textdef{transverse} (or \textdef{\( 1 \)-transverse}) if it satisfies the condition~\( RI \cap S = I \).
	For~\( n \ge 2 \), it is called \textdef{\( n \)-transverse} if it is \( (n - 1) \)-transverse
	and furthermore satisfies~\( RI^{n} \cap S = I^{n} \).
\end{definition}

We shall later see that our particular choice of involution~\( \omega \) is both transverse and \( 2 \)-transverse.
We leave the questions of higher transversality as a conjecture.

\begin{proposition}\label{res:fundamental_ses}
	For a transverse involution~\( \omega \),
	we have a short exact sequence of vector spaces
	\[\begin{tikzcd}[sep=small]
		0 \ar[r]
			&
			S/I
			\ar[r]
				&
				R/RI
				\ar[r, "\diff{d}"]
					& 
					I/I^{2}
					\ar[r]
						& 0
	\text{.}
	\end{tikzcd}\]
\end{proposition}

\begin{proof}
	\zcref[S]{res:preimage_of_I^(n+1)} below shows that
	the kernel of~\(
		\diff{d}
		\colon
		R/RI
		\to
		I/I^{2}
	\) is
	the image of~\( S \)
	under the quotient map~\(
		R
		\to
		R/RI
	\),
	and that image is~\(
		S/(RI \cap S)
	\).
	But this is equal to~\( S/I \) due to transversality of~\( \omega \).
\end{proof}

\begin{lemma}\label{res:preimage_of_I^(n+1)}
	The preimage of~\( I^{n +1} \) under~\( \diff{d} \)
	is~\( S + RI^{n} \) for all~\( n \ge 1 \).
\end{lemma}

\begin{proof}
	If~\( \diff{d\differentialdspace} x \in I^{n +1} \), we may write~\(
		\diff{d\differentialdspace} x 		=
		\sum y_{i} \diff{d\differentialdspace} z_{i} 	\)
	for~\( y_{i} \in I^{n} \) and~\( z_{i} \in R \).
	Due to the \( S \)-linearity of~\( \diff{d} \),
	this means that~\(
		\diff{d\differentialdspace} x 		=
		\diff{d} (\sum y_{i} z_{i} )
	\),
	i.e.~\(
		\diff{d} (x +\sum y_{i} z_{i} ) = 0
	\).
	Therefore,~\(
		x
		+
		\sum y_{i} z_{i} 		\in
		S
	\),
	and consequently,~\(
		x \in S + RI^{n}
	\).
\end{proof}

\begin{notation}
	Given a \( \mathbb{Z}_{\ge 0} \)-graded ring~\smash{\(
		R
		=
		\bigoplus_{i \ge 0} R_{i} 	\)},
	we denote by~\( R_{+} \)
	the augmentation ideal~\smash{\(
		R_{+}
		=
		\bigoplus_{i \ge 1} R_{i} 	\)}.
\end{notation}

\begin{definitionbreak}\label{def:normal}
	\begin{definitionlist}
		\item The \textdef{norm} of the involution~\( \omega \)
		is defined by~\( N(x) = x \, \widebar{x} \) for~\( x \in R \).
		Clearly, we have~\( N(x) \in S \) for all~\( x \).
		\item Suppose that \( R \)~is a \( \mathbb{Z}_{\ge 0} \)-graded ring.
		If we have~\( N(x) \in I \) for all~\( x \in R_{+} \),
		we say that the involution~\( \omega \)
		is \textdef{normal}.
	\end{definitionlist}
\end{definitionbreak}

Again, we shall see that this is the case for our particular choice of~\( \omega \).

\begin{proposition}\label{res:normality_enough_for_generators}
	The norm map defines a map of rings~\( N \colon R \to S/I \).
	Therefore, if \( R \)~is \( \mathbb{Z}_{\ge 0} \)-graded,
	normality of~\( \omega \)~can be proved by checking the condition~\( N(x) \in I \)
	on a set of generators for~\( R_{+} \).
\end{proposition}

\begin{proof}
	The Thick Leibniz Rule~\(
		\diff{d} (x y )
		=
		x \diff{d\differentialdspace} y 		+
		y \diff{d\differentialdspace} x 		+
		\diff{d\differentialdspace} x \diff{d\differentialdspace} y 	\)
	shows that~\(
		x \diff{d\differentialdspace} y 		+
		y \diff{d\differentialdspace} x 		\in
		I
	\)
	for all~\( x, y \in R \).
	This implies that~\( N(x +y) = N(x) + N(y) \)
	in~\( S/I\).
	All other conditions for~\( N \)
	to be a ring map are obvious.
\end{proof}

\begin{definition}
	For a \( \mathbb{Z}_{\ge 0} \)-graded ring~\( R \),
	the \textdef{ideal of squares}~\( Q(R) \)
	is the ideal generated by~\( x^{2} \) for all~\( x \in R_{+} \).
\end{definition}

\begin{proposition}\label{res:z^2_in_RI}
	If \( R \)~is \( \mathbb{Z}_{\ge 0} \)-graded and the involution~\( \omega \) is normal,
	then we have~\( x^{2} \in RI \)
	for all~\( x \in R_{+} \).
	Therefore,~\(
		Q(R)
		\subset
		RI
	\).
\end{proposition}

\begin{proof}
	Follows from~\(
		x^{2}
		=
		x \, \widebar{x} 		+
		x \diff{d\differentialdspace} x 	\).
\end{proof}

\chapter{Basic properties of the ring~\texorpdfstring{\( R \)}{R}}

The main character of this story will be the
polynomial ring in infinitely many variables
\[
	R = \mathbb{F}_{2} \lbrack w_{1}, w_{2}, w_{3}, \dotsc \rbrack
\]
endowed with a \( \mathbb{Z}_{\ge 0} \)-grading
given by~\( \deg (w_{i}) = i \).
It is equipped with an involution~\(
	\omega
	\colon
	R
	\to
	R
\)
which we write as~\( x \mapsto \widebar{x} \),
and which is given recursively by the formula
\begin{equation}\label{eq:omega_definition}
	\sum_{i =0}^{k} w_{i} \widebar{w}_{k -i} 	=
	0
\end{equation}
where we use the convention~\( \widebar{w}_{0} = w_{0} = 1 \).
We define the map~\( \diff{d} \colon R \to R \),
the subring~\( S \subset R \), and the ideal~\( I \subset S \) as in \zcref{chap:general_ring_vocabulary}.

Applying the formula~\eqref{eq:omega_definition} inductively, we arrive at the formula
\begin{equation}\label{eq:omega_sum_formula}
	\widebar{w}_{k}
	=
	\sum_{i_{1} +\dotsb +i_{p} =k} w_{i_{1}} \dotsm w_{i_{p}} 	\text{.}
\end{equation}
This can be conveniently rewritten by introducing the formal sum
\[
	W = w_{1} +w_{2} +w_{3} +\dotsb 	\text{.}
\]
Then the formula~\eqref{eq:omega_sum_formula} simply says that
\begin{equation}\label{eq:formal_omega_w}
	\widebar{W}
	=
	\sum_{r =1}^{\infty} W^{r} 	=
	W \sum_{r =0}^{\infty} W^{r} 	=
	W (1+\widebar{W}) 	\text{,}
\end{equation}
which implies the formulae
\begin{equation}\label{eq:formal_w_sum=product}
	W + \widebar{W}
	=
	W \, \widebar{W} 	=
	W^{2} (1+\widebar{W}) 	\text{.}
\end{equation}
If we extend the formal sum as~\(
	W_{\!+}
	=
	1+w_{1} +w_{2} +w_{3} +\dotsb \),
we also arrive at the compact expression
\begin{equation}\label{eq:formal_whit_inverse}
	W_{\!+} \, \widebar{W}_{\!\!\!+} 	=
	1
	\text{.}
\end{equation}

We are now ready to prove that the involution~\( \omega \) is normal
in the sense of \zcref{def:normal}:

\begin{proposition}\label{res:omega_is_normal}
	The involution~\( \omega \)
	on~\( R \) is normal.
\end{proposition}

\begin{proof}
	\zcref[S]{res:normality_enough_for_generators} shows that it is enough
	to check that~\( N(w_{j}) \in I \)
	for all~\( j \).
	But this follows from~\eqref{eq:omega_definition}
	with~\( k = 2j \)
	since~\(
		w_{i} \widebar{w}_{k -i} 		+
		w_{k -i} \widebar{w}_{i} 		=
		\diff{d} (w_{i} \widebar{w}_{k -i} )
	\)
	for all~\( i \neq j \).
\end{proof}

\begin{proposition}\label{res:Q(R)=RI}
	The ideal of squares of~\( R \)
	is~\( Q(R) = RI \).
\end{proposition}

\begin{proof}
	We already have the inclusion~\(
		Q(R)
		\subset
		RI
	\)
	from \zcref{res:z^2_in_RI}.
	To see the converse,
	note that it follows from the Thick Leibniz Rule
	that \( RI \)~is generated as an \( R \)-module
	by the elements~\( \diff{d\differentialdspace} w_{i} \).
	It is therefore enough to prove that these
	lie in~\( Q(R) \).
	But this follows by writing the formula~\eqref{eq:formal_w_sum=product}
	out in degree~\( i \).
\end{proof}

\begin{remark}\label{rem:dw_2i_related_to_w_i^2}
	For later use, let us explicitly write out~\eqref{eq:formal_w_sum=product}
	in degree~\( 2i \):
	\[
		\diff{d\differentialdspace} w_{2i } 		=
		w_{i}^{2}
		+
		\sum_{k =1}^{i -1} \widebar{w}_{2(i -k) } \, w_{k}^{2} 		\text{.}
	\]
\end{remark}

\zcref{res:Q(R)=RI} shows that squares play an important role in our study of \( I \) and~\( RI \).
In the following sections, square-free monomials will appear regularly.

\chapter{The subring of power sums}

An important class of elements of~\( S \) is the \textdef{power sums},
denoted by~\( p_{k} \) for~\( k \ge 1 \).
Inside the rings~\( R_{(n)} \) from
the introduction, these are defined as~\(
	p_{k}
	=
	\smash{\sum x_{i}^{k} }
\).
In~\( R \), they can, for instance, be defined recursively using the Newton identity
\parencite[see e.g.][Exercise~1, section~6.1]{ful}
\begin{equation}\label{eq:newton_identity}
	k w_{k} 	=
	\sum_{j =0}^{k -1} p_{k -j} w_{j} 	\mathrlap{\qquad\qquad(k \ge 1)\text{.}}
\end{equation}
Due to characteristic~\( 2 \), they satisfy the additional property
\begin{equation}\label{eq:p_i^2=p_(2i)}
	p_{k}^{2}
	=
	p_{2k }
	\text{,}
\end{equation}
which holds because it holds in all~\( R_{(n)} \).
We denote by~\( \mathit{PS} \subset S \)
the subring generated by the power sums~\( p_{k} \).

\begin{proposition}\label{res:odd_and_even_formulae}
	We have the following formulae for~\( p_{k} \) for all~\( k \ge 1 \):
	\[
		p_{k}
		=
		\sum_{\substack{0\le i \le k \\i \textup{\space odd}}} w_{i} \widebar{w}_{k -i} 		=
		\sum_{\substack{0\le i \le k \\i \textup{\space even}}} w_{i} \widebar{w}_{k -i} 		\text{.}
	\]
\end{proposition}

\begin{proof}
	We introduce the additional formal sums \(
		W_{\mathrm{odd}}
		=
		w_{1} +w_{3} +w_{5} +\dotsb 	\)
	and~\(
		P
		=
		p_{1} +p_{2} +p_{3} +\dotsb 	\).
	Then the Newton identity~\eqref{eq:newton_identity}
	can be rewritten as
	\[
		P \, W_{\!+} 		=
		W_{\mathrm{odd}}
		\text{.}
	\]
	Because of~\eqref{eq:formal_whit_inverse},
	this is equivalent to~\(
		P
		=
		W_{\mathrm{odd}} \, \widebar{W}_{\!\!\!+} 	\),
	which, when written out degreewise, yields the first formula.
	The second one then follows from the first one using~\eqref{eq:omega_definition}.
	(Alternatively, one can give a proof in terms of Young diagrams using~\eqref{eq:diagram_formula_for_p_k} below. While less elegant, this approach has the benefit of proving the formulae in arbitrary characteristic.)
\end{proof}

\begin{remark}
	This provides a generalization and explanation for
	Theorem~1.2 (and Theorem~4.6)
	in \textcite{matsz_wendt}.
\end{remark}

\begin{corollary}\label{res:p^2_formula}
	For odd~\( k \ge 1 \), we have the formula
	\[
		p_{k}^{2}
		=
		\sum_{i =0}^{(k -1)/2} \bigl[\diff{d} (w_{2i } w_{2(k -i) } ) +\diff{d\differentialdspace} w_{2i } \diff{d\differentialdspace} w_{2(k -i) } \bigr] 		\text{.}
	\]
\end{corollary}

\begin{proof}
	By applying the \enquote{even} formula
	from \zcref{res:odd_and_even_formulae} to~\( p_{k}^{2} = p_{2k } \), we get
	\[
		p_{k}^{2}
		=
		\sum_{i =0}^{(k -1)/2} \diff{d} (w_{2i } \widebar{w}_{2(k -i) }) 		\text{.}
	\]
	Then the result follows by plugging in~\(
		\widebar{w}_{2(k -i) }
		=
		w_{2(k -i) }
		+
		\diff{d\differentialdspace} w_{2(k -i) } 	\).
\end{proof}

\begin{corollary}\label{res:formula_for_w^2}
	For any~\( k \ge 1 \),
	we have
	\[
		w_{2k }^{2}
		=
		w_{2k } \diff{d\differentialdspace} w_{2k } 		+
		p_{k}^{4}
		+
		\sum_{i =0}^{k -1} \;\;\bigl[\diff{d} (w_{2i } w_{4k -2i } ) +\diff{d\differentialdspace} w_{2i } \diff{d\differentialdspace} w_{4k -2i } \bigr] 		\text{.}
	\]
\end{corollary}

\begin{proof}
	Applying the \enquote{even} variant of~\zcref{res:p^2_formula} to~\( p_{k}^{4} = p_{4k } \),
	we get
	\[
		p_{k}^{4}
		=
		p_{4k }
		=
		w_{2k } \widebar{w}_{2k } 		+
		\sum_{i =0}^{k -1} \diff{d} (w_{2i } \widebar{w}_{4k -2i } ) 		\text{.}
	\]
	Then the result follows by plugging in~\(
		w_{2k } \widebar{w}_{2k } 		=
		w_{2k }^{2}
		+
		w_{2k } \diff{d\differentialdspace} w_{2k } 	\)
	and~\(
		\widebar{w}_{4k -2i }
		=
		w_{4k -2i }
		+
		\diff{d\differentialdspace} w_{4k -2i } 	\).
\end{proof}

\chapter{The quotient ring~\texorpdfstring{\( S/I \)}{S/I}}

As mentioned in the introduction,
the classical theory of symmetric functions
provides us with an \( \mathbb{F}_{2} \)-basis for~\( R \) consisting of the
Schur functions~\(
	s_{\lambda}
\)
corresponding to partitions~\( \lambda \).
In these terms, the involution~\( \omega \)
is given by~\(
	\widebar{s}_{\lambda}
	=
	s_{\lambda^{\!\vee}}
\),
where \( \lambda^{\!\vee} \)~denotes the conjugate partition
whose Young diagram is obtained by mirroring the Young diagram
of~\( \lambda \).
Therefore, the subring~\( S = R^{\omega} \subset R \)
of \( \omega \)-invariants
has a basis consisting of
of~\(
	s_{\lambda}
\)
for~\( \lambda = \lambda^{\!\vee} \)
as well as~\(
	\diff{d\differentialdspace} s_{\lambda} 	=
	s_{\lambda}
	+
	s_{\lambda^{\!\vee}}
\)
for~\( \lambda \neq \lambda^{\!\vee} \).
Hence the Schur functions~\( s_{\lambda} \),
for self-conjugate partitions~\( \lambda \),
form a basis for the quotient ring~\( S/I \).

One notable class of self-conjugate partitions are the following, which will later serve as building blocks for all self-conjugate partitions:
\begin{definition}
	Given a number~\( j \), we denote by~\( \Gamma (j) \)
	the partition whose Young diagram is a symmetric hook with \( 2j - 1 \)~boxes, i.e.
	\[
		\Gamma (j)
		\quad
		=
		\quad
		\begin{ytableau}
		{}
              & \none[\cdots]
                  &
		\\
          \none[\raisebox{-.2em}{\vdots}]
		\\
          {}
		\end{ytableau}
	\]
	where the side length is~\( j \).
\end{definition}

The \ABBRVMNrule{}
\parencite[Theorem~7.17.1]{stanley}
says that we can calculate the product~\( p_{k} \, s_{\lambda} \) as
\begin{equation}\label{eq:murnaghan-nakayama_rule}
	\textstyle
	p_{k} \, s_{\lambda} 	=
	\sum_{\mu} s_{\mu} 	\text{,}
\end{equation}
where the sum runs over all partitions~\( \mu \)
containing~\( \lambda \) and
such that the diagram quotient~\( \mu /\lambda \)
(which is the skew diagram obtained by erasing~\( \lambda \)
from~\( \mu \))
is a border strip diagram of size~\( k \).
In particular, taking~\( \lambda \) to be the
empty diagram, we get that
\begin{equation}\label{eq:diagram_formula_for_p_k}
	\textstyle
	p_{k}
	=
	\sum_{\mu} s_{\mu} \end{equation}
with~\( \mu \) running
over all hook-shaped Young diagrams of size~\( k \).

\begin{lemma}\label{res:p_k_in_S/I}
	We have
	\(
		p_{2j -1}
		=
		s_{\Gamma (j)}
		\in S/I
	\)
	and~\(
		p_{2j }
		=
		p_{j}^{2}
		=
		0
		\in
		S/I
	\) for all~\( j \).
\end{lemma}

\begin{proof}
	The formula~\( p_{2j } = p_{j}^{2} \)
	comes from~\eqref{eq:p_i^2=p_(2i)},
	and the rest follows from~\eqref{eq:diagram_formula_for_p_k}
	since we have~\(
		s_{\mu}
		+
		s_{\mu^{\vee}}
		=
		\diff{d\differentialdspace} s_{\mu} 		\in
		I
	\).
\end{proof}

\begin{lemma}\label{res:p1...pk=square}
	We have~\(
		p_{1} p_{3} \dotsm p_{2j-1} 		=
		s_{\operatorname{Sq} (j)}
		\in
		S/I
	\)
	where \( \operatorname{Sq} (j) \)~is the partition
	whose Young diagram is the \( j \times j \)~square.
\end{lemma}

\begin{proof}
	We claim that~\(
		p_{2j-1} \, s_{\operatorname{Sq} (j -1)} 		=
		s_{\operatorname{Sq} (j)}
		\in
		S/I
	\).
	Indeed, the \ABBRVMNrule{}~\eqref{eq:murnaghan-nakayama_rule}
	shows that it is the sum of~\( s_{\mu} \)
	for all self-conjugate~\( \mu \)
	containing~\( \operatorname{Sq} (j -1) \)
	for which \( \mu /\operatorname{Sq} (j -1) \)~is
	a border strip diagram of size~\( 2j - 1 \).
	But the only self-conjugate diagram we can get by adding such a border strip
	to~\(
		\operatorname{Sq} (j -1)
	\)
	is~\(
		\mu
		=
		\operatorname{Sq} (j)
	\).
\end{proof}

This allows us to give a description of the product structure in~\( S/I \):

\begin{proposition}\label{res:S/I_exterior_algebra}
	The ring~\( S/I \)
	is the exterior algebra over~\( \mathbb{F}_{2} \)
	on the vector space spanned by~\( p_{k} \) for~\( k \)~odd.
	In other words,
	the map~\(
		\mathit{PS}/Q(\mathit{PS})
		\to
		S/I
	\)
	is an isomorphism.
\end{proposition}

\begin{proof}
	Due to \zcref{res:p_k_in_S/I}, we have~\( p_{k}^{2} = 0 \in S/I \),
	which shows that there is a map~\(
		\bigwedge (p_{1}, p_{3}, p_{5}, \dotsc)
		\to
		S/I
	\).
	For odd~\( k \),
	we get from
	\zcref{res:p1...pk=square} that the map~\(
		\bigwedge (p_{1}, p_{3}, \dotsc, p_{k})
		\to
		S/I
	\)
	is non-zero on the top form.
	By \zcref{res:all_divide_top_form} below,
	all non-zero elements divide the top form,
	so this map must be injective.
	Therefore, \(
		\bigwedge (p_{1}, p_{3}, p_{5}, \dotsc)
		\to
		S/I
	\)~is injective as well.
	
	Because of the description of the product of Schur functions in terms
	of \ABBRVlrcoef{}s, the image of~\(
		\bigwedge (p_{1}, p_{3}, \dotsc, p_{k})
		\to
		S/I
	\)
	is a subspace of the vector space spanned by~\( s_{\lambda} \)
	for self-conjugate~\( \lambda \) contained inside~\(
		\operatorname{Sq} ((k +1)/2)
	\).
	By comparing dimensions, we see that the image must be the whole of this space. Therefore, \(
		\bigwedge (p_{1}, p_{3}, p_{5}, \dotsc)
		\to
		S/I
	\)~is an isomorphism.
\end{proof}

\begin{lemma}[Linear-Algebraic Lemma]\label{res:all_divide_top_form}
	Let~\( V \) be a vector space of dimension~\( n < \infty \)
	over the field~\( \mathbb{F}_{2} \),
	and let~\( e \in \bigwedge^{\!\!n} (V) \) be the non-zero top form.
	Then any non-zero~\( x \in \bigwedge (V) \)
	divides~\( e \).
\end{lemma}

This is, of course, a standard linear-algebraic fact, but for some reason,
I could not find it in the literature except in characteristic~\( 0 \).

\begin{proof}
	We choose a basis~\( v_{1}, \dotsc, v_{n} \) for~\( V \)
	and
	write~\( x \)
	as a linear combination of the wedge basis vectors~\(
		v_{i_{1}} \wedge \dotsb \wedge v_{i_{k}} 	\)
	with~\( i_{1}<\dotsb <i_{k} \).
	Let~\(
		v_{i_{1}} \wedge \dotsb \wedge v_{i_{k}} 	\)
	be a basis vector from this decomposition with the additional property
	that it is not divisible by any of the other vectors
	from the decomposition.
	We write~\(
		j_{1}<\dotsb <j_{n -k} 	\)
	for the indices that did not occur in the sequence~\(
		i_{1}<\dotsb <i_{k} 	\)
	and put~\(
		y
		=
		v_{j_{1}} \wedge \dotsb \wedge v_{j_{n -k}} 	\)
	Then~\(
		(v_{i_{1}} \wedge \dotsb \wedge v_{i_{k}} ) \wedge y 		= e
	\)
	while any other basis vector in the decomposition is annihilated
	by~\( y \).
	Therefore,~\(
		x \wedge y 		=
		e
	\).
\end{proof}

\chapter{Mixed coordinates, the Transversality Theorems, and standard forms}\label{chap:mixed_coordinates}

It turns out that we get an equally good set of coordinates for~\( R \) by
replacing~\( w_{k} \) for odd~\( k \) by~\( p_{k} \).

\begin{proposition}\label{res:mixed_coordinates}
	The ring~\( R \) is polynomially generated by \( p_{2i -1} \) and~\( w_{2i } \) for~\( i \ge 1 \).
	In other words, \( R \)~can be regarded as a polynomial algebra
	\[
		R
		=
		\mathbb{F}_{2} \lbrack p_{1}, w_{2}, p_{3}, w_{4}, p_{5}, w_{6}, \dotsc \rbrack
		=
		\mathit{PS} \lbrack w_{2}, w_{4}, w_{6}, \dotsc \rbrack
	\text{.}
	\]
\end{proposition}

\begin{definition}\label{def:mixed_coordinates}
	We shall refer to these coordinates as \textdef{mixed coordinates}.
\end{definition}

\begin{proof}
	Using the Newton identity~\eqref{eq:newton_identity}
	inductively together with~\eqref{eq:p_i^2=p_(2i)}, we get that \( w_{k} \) for odd~\( k \)
	can be expressed in the alternative coordinates.
	This shows that the map~\(
		\mathit{PS} \otimes_{\mathbb{F}_{2}} \mathbb{F}_{2} \lbrack w_{2}, w_{4}, w_{6}, \dotsc \rbrack 		\to
		R
	\)
	is surjective.
	Since both the source and target
	are polynomial algebras with
	one generator in each positive degree, the dimensions agree in each degree,
	so the map must in fact be an isomorphism.
\end{proof}

This allows us to prove:

\begin{theorem}[(First) Transversality Theorem]\label{res:transversality_theorem}
	The involution~\( \omega \) on~\( R \)
	is transverse.
\end{theorem}

In a previous version of this article, this theorem was essential to the proof of \zcref{res:generatorss_and_relations_stable_case}.
However, after simplifying that proof, it is no longer needed there. We nevertheless believe that the Transversality Theorem is of independent interest and have therefore chosen to keep it.

\begin{proof}
	\zcref[S]{res:mixed_coordinates}
	shows that
	the map
	\[
		\mathit{PS}/Q(\mathit{PS})
		\to
		R/Q(R)
	\]
	is an inclusion.
	Furthermore,
	\zcref{res:S/I_exterior_algebra,res:Q(R)=RI}
	tell us that
	\(
		R/Q(R)
		=
		R/RI
	\)
	and~\(
		\mathit{PS}/Q(\mathit{PS})
		=
		S/I
	\),
	respectively.
	Hence this inclusion factors as
	\[
		S/I
		\to
		S/(RI \cap S)
		\subset
		R/RI
		\text{.}
	\]
	Therefore, the inclusion~\(
		I
		\subset
		RI \cap S
	\)
	must be an equality.
\end{proof}

In analogy with differential forms, we introduce the following notation:

\begin{definition}
	We write~\( (\Omega^{\SemantexBullet}, \diff{d}) \) for the complex
	given by~\(
		\Omega^{n} = RI^{n}/RI^{n +1}
	\)
	for all~\( n \ge 0 \)
	and with the differential~\(
		\diff{d}
		\colon
		\Omega^{n}
		\to
		\Omega^{n +1}
	\)
	induced by the map~\( \diff{d} \).
\end{definition}

\begin{proposition}\label{res:dw_2i_basis_for_RI^n/RI^(n+1)}
	For any~\( n \),
	\( \Omega^{n} \)~is
	a free \( R/RI \)-module with
	basis represented by the products~\(
		\diff{d\differentialdspace} w_{2i_{1}} \dotsm \diff{d\differentialdspace} w_{2i_{n}} 	\)
	with~\(
		i_{1} \le \dotsb \le i_{n} 	\).
\end{proposition}

\begin{proof}
	Using the description~\(
		RI^{n}/RI^{n +1}
		=
		Q(R)^{n}/Q(R)^{n +1}
	\) (c.f.~\zcref{res:Q(R)=RI}),
	we note that we have a basis represented by~\smash{\(
		w_{i_{1}}^{2} \dotsm w_{i_{n}}^{2} 	\)}
	for~\(
		i_{1} \le \dotsb \le i_{n} 	\).
	Also, \zcref{rem:dw_2i_related_to_w_i^2}
	shows that \( \diff{d\differentialdspace} w_{2i } \)~is~\( w_{i}^{2} \)
	plus a linear combination of~\( w_{k}^{2} \)
	for~\( k < i \).
	Therefore, a grading argument shows that~\(
		\diff{d\differentialdspace} w_{2i_{1} } \dotsm \diff{d\differentialdspace} w_{2i_{n} } 	\)
	for~\(
		i_{1} \le \dotsb \le i_{n} 	\)
	also form a basis.
\end{proof}

The analogy with differential forms is justified by the following observation,
which will not be used any further in this paper:

\begin{proposition}
	We have~\smash{\(
		\Omega^{1}
		\cong
		\Omega_{\mathit{PS}}^{1} (R/RI)
	\)},
	where the right-hand side is understood as the
	module of Kähler differential \( 1 \)-forms as defined
	in e.g.~\textcite[section~II.8]{hartshorne}.
\end{proposition}

\begin{proof}
	Proposition~II.8.4A in \textcite{hartshorne}
	provides us with an exact sequence of \( R/RI \)-modules
	\[\begin{tikzcd}[sep=small]
		RI/RI^{2}
		\ar[r, "\delta"]
			&
			\Omega_{\mathit{PS}}^{1} (R) \otimes_{R} R/RI 			\ar[r]
				& \Omega_{\mathit{PS}}^{1} (R/RI)
				\ar[r]
					& 0
	\end{tikzcd}\]
	where~\( \delta (x) = \diff{d\differentialdspace} x \otimes 1\).
	Also,
	\zcref{res:mixed_coordinates} shows that \( R \)~is a polynomial algebra over~\( \mathit{PS} \).
	Hence
	Example~II.8.2.1 in \textcite{hartshorne}
	tells us that
	\smash{\( \Omega_{\mathit{PS}}^{1} (R) \)}~is the free
	\( R \)-module generated by~\( \diff{d\differentialdspace} w_{2i } \)
	for all~\( i \).
	Therefore, \smash{\(
		\Omega_{\mathit{PS}}^{1} (R) \otimes_{R} R/RI 	\)}~is a free \( R/RI \)-module with the same basis.
	Comparing this to \zcref{res:dw_2i_basis_for_RI^n/RI^(n+1)},
	we get
	\smash{\(
		\Omega_{\mathit{PS}}^{1} (R) \otimes_{R} R/RI 		\cong
		\Omega^{1}
	\)}.
	In other words, we are done if we can prove that the map~\( \delta \) is zero.
	Because of \zcref{res:Q(R)=RI},
	this amounts to proving that~\( \delta (a x^{2}) = 0 \)
	for all~\( a, x \in R \).
	Since~\( \diff{d} (x^{2}) = 2 x \diff{d\differentialdspace} x = 0 \),
	we have~\( \diff{d} (a x^{2}) = x^{2} \diff{d\differentialdspace} a \).
	But~\smash{\(
		\delta (a x^{2} )
		=
		x^{2} \diff{d\differentialdspace} a \otimes 1		=
		0
		\in
		\Omega_{\mathit{PS}}^{1} (R) \otimes_{R} R/RI 	\)}
	since~\( x^{2} \in Q(R) = RI \).
\end{proof}

\begin{lemma}\label{res:kernel_of_d_on_Omega^0}
	The kernel of~\(
		\diff{d}
		\colon
		\Omega^{0} = R/RI
		\to
		\Omega^{1} = RI/RI^{2}
	\)
	is~\( S/I \).
\end{lemma}

Here, we regard~\( S/I \) as a subring of~\( R/RI \)
because of \zcref{res:fundamental_ses}.

\begin{proof}
	Due to \zcref{res:Q(R)=RI,res:mixed_coordinates},
	we can represent~\( x \in R/RI \) as a polynomial in mixed coordinates where all of the monomials are square-free.
	Note that the map~\(
		\diff{d}
		\colon
		\Omega^{0}
		\to
		\Omega^{1}
	\)
	is a derivation in the classical sense
	since the term~\( \diff{d\differentialdspace} x \diff{d\differentialdspace} y \)
	in the formula~\eqref{eq:thick_leibniz} vanishes.
	Therefore, because of \zcref{res:dw_2i_basis_for_RI^n/RI^(n+1)}, we get the formula
	\[
		\diff{d\differentialdspace} x 		=
		\sum \frac{\partial x}{\partial w_{2i }} \diff{d\differentialdspace} w_{2i } 		\text{,}
	\]
	which shows that~\( \diff{d\differentialdspace} x = 0 \) if and only if~\(
		\partial x /\partial w_{2i }
		\in
		RI
		=
		Q(R)
	\) for all~\( i \).
	For a polynomial without squares in characteristic~\( 2 \),
	this must mean that \( x \)~was in fact only a polynomial in the variables~\( p_{2i -1} \).
	So~\( x \in S/I \), and we are done.
\end{proof}

We include the following for completeness even though it will not be needed any further in this paper:

\begin{theorem}[Second Transversality Theorem]
	The involution~\( \omega \)
	is \( 2 \)-transverse.
\end{theorem}

\begin{proof}
	Tautologically, the image of~\(
		\diff{d}
		\colon
		\Omega^{0}
		\to
		\Omega^{1}
	\)
	is~\(
		I/(RI^{2} \cap I)
		\subset
		RI/RI^{2}
	\).
	Furthermore,
	\( \diff{d} \)
	factors as~\(
		R/RI
		\to
		I/I^{2}
		\to
		I/(RI^{2} \cap I)
	\).
	\zcref[S]{res:fundamental_ses,res:kernel_of_d_on_Omega^0}
	show that 
	\(
		R/RI
		\to
		I/I^{2}
	\)~has the same kernel as~\(
		R/RI
		\to
		I/(RI^{2} \cap I)
	\),
	hence the map~\(
		I/I^{2}
		\to
		I/(RI^{2} \cap I)
	\)
	is an equality.
	Therefore,~\(
		I^{2}
		=
		RI^{2} \cap I
	\),
	and because we already know
	from the First Transversality Theorem
	that~\( I = RI \cap S \),
	this proves that~\(
		I^{2}
		=
		RI^{2} \cap S
	\).
\end{proof}

We expect that higher transversality also holds and leave this as a conjecture for now:

\begin{conjecture}[Conjecture]
	The involution~\( \omega \)
	is also \( n \)-transverse for~\( n \ge 3 \).
\end{conjecture}

Finally, we use \zcref{res:dw_2i_basis_for_RI^n/RI^(n+1)}
to construct a basis for all of~\( R \)
which will come in handy when proving the main theorem of this paper:

\begin{theorem}\label{res:existence_of_the_standard_form}
	We have an \( \mathbb{F}_{2} \)-basis for~\( R \)
	consisting of all expressions of the form
	\[
		p a \diff{d\differentialdspace} w_{2i_{1} } \dotsm \diff{d\differentialdspace} w_{2i_{n} } 		\text{,}
	\]
	where~\( p \) is a non-zero, square-free monomial
	in the variables~\( p_{2i -1} \), 
	while \( a \)~is a non-zero, square-free monomial in the variables~\( w_{2i } \),
	and~\(
		i_{1} \le \dotsb \le i_{n} 	\)
	with~\( n \ge 0 \).
\end{theorem}

\begin{definition}
	When an element is written in this basis, we call the resulting
	expression the \textdef{standard form} of that element.
\end{definition}

\begin{proof}
    To prove that the system is linearly independent,
    we let~\( \sum x_{i} = 0 \) be a non-trivial linear relation between them.
    Let~\( n \) be the maximal number such that~\(
        x_{i} \in RI^{n}
    \)
    for all~\( i \).
    By assumption, the set~\( J \)
    of all~\( i \) satisfying~\(
        x_{i} \neq 0 \in \Omega^{n}
    \)
    is non-empty.
    Therefore, we find that the linear relation~\( \sum x_{i} = 0 \in R \) descends to the linear relation~\( \sum_{i \in J} x_{i} = 0 \in \Omega^{n} \).
    \zcref[S]{res:dw_2i_basis_for_RI^n/RI^(n+1)}
    shows that~\(
    	\diff{d\differentialdspace} w_{2i_{1} } \dotsm \diff{d\differentialdspace} w_{2i_{n} }     \)
    with~\(
    	i_{1} \le \dotsb \le i_{n}     \) constitute a basis
    for~\( \Omega^{n} \) over~\( R/RI \).
    Therefore, the same space, regarded as a vector space over~\( \mathbb{F}_{2} \),
    has a basis consisting of~\(
        p a \diff{d\differentialdspace} w_{2i_{1} } \dotsm \diff{d\differentialdspace} w_{2i_{n} }     \),
    where \( p \) and~\( a \) run over all non-zero, square-free monomials in the variables \( p_{2i -1} \) and~\( w_{2i } \), respectively.
    In conclusion, the linear relation~\( \sum_{i \in J} x_{i} = 0 \) must be trivial, hence~\( J = \varnothing \), and we have a contradiction.

    To prove that they span all of~\( R \),
    we note that we clearly do get a basis if we replace~\( \diff{d\differentialdspace} w_{2i_{j} } \) by~\smash{\( w_{i_{j}}^{2} \)}
    in the definition, i.e.\ look
    at~\smash{\(
    	p a \, w_{i_{1}}^{2} \dotsm w_{i_{n}}^{2}     \)}
    instead.
    Therefore, it suffices to prove that any such element is contained in the span.
    But this follows from \zcref{rem:dw_2i_related_to_w_i^2}
    using a grading argument.
\end{proof}

\chapter{Generators and relations for~\texorpdfstring{\( S \)}{S}}\label{chap:generators_and_relations_for_S}

In this section, we shall provide a presentation by generators and relations for~\( S \) as an algebra over~\( \mathit{PS} \).
The starting point is the observation that applying~\( \diff{d} \)
to the Thick Leibniz Rule as stated in~\eqref{eq:thick_leibniz_general} yields the following relation in~\( S \):
\begin{equation}\label{eq:d^2_thick_Leibniz_relation_inspiration}
    0
    =
    \diff{d}^{2} (x_{1} \dotsm x_{n} )
    =
    \sum_{\varnothing \neq T \subsetneq \lbrace 1, \dotsc, n \rbrace}
    \diff{d} (x_{T^{\complement}}) \, (\diff{d\differentialdspace} x )_{T}     \text{.}
\end{equation}
We shall see that these exhaust most of the relations in~\( S \).
In fact, only a small subset of them will be necessary to generate the others.

We denote by~\( F \)
the free commutative \( \mathit{PS} \)-algebra generated
by symbols~\( \diff{\delta \differentialdeltaspace} a \)
where \( a \)~is a square-free monomial of strictly positive degree in the variables~\( w_{2i } \).
In other words,
we only use~\(
	a
	=
	w_{2i_{1} } \dotsm w_{2i_{n} } \)
where~\(
	i_{1}<\dotsb <i_{n} \)
and~\( n \ge 1 \).

We then extend the notation \( \diff{\delta \differentialdeltaspace} x \)
to arbitrary~\( x \in R \)
using \zcref{res:existence_of_the_standard_form}:
If~\(
	x
	=
	p a \diff{d\differentialdspace} w_{2i_{1} } \dotsm \diff{d\differentialdspace} w_{2i_{n} } \),
where \( p \)~is a square-free monomial in the variables~\( p_{2i -1} \),
while \( a \)~is a square-free monomial in the variables~\( w_{2i } \),
and~\(
	i_{1} \le \dotsb \le i_{n} \),
we then define~\(
	\diff{\delta \differentialdeltaspace} x 	=
	p \diff{\delta \differentialdeltaspace} a \diff{\delta \differentialdeltaspace} w_{2i_{1} } \dotsm \diff{\delta \differentialdeltaspace} w_{2i_{n} } \).
In the case~\( a = 1 \), we use the understanding
that~\( \diff{\delta \differentialdeltaspace} a = 0 \), and hence~\( \diff{\delta \differentialdeltaspace} x = 0 \).
Then we extend this linearly to arbitrary~\( x \in R \).

For the relations,
we take inspiration from~\eqref{eq:d^2_thick_Leibniz_relation_inspiration}
and define an operation~\( \delta^{2} \):
Given a collection~\( (x_{1}, \dotsc, x_{n}) \)
of elements of~\( R \), we define~\( \delta^{2} (x_{1}, \dotsc, x_{n}) \in F \) by the formula
\[
        \delta^{2} (x_{1}, \dotsc, x_{n})
        =
        \sum_{\varnothing \neq T \subsetneq \lbrace 1, \dotsc, n \rbrace}
        \diff{\delta} (x_{T^{\complement}}) \, (\diff{\delta \differentialdeltaspace} x )_{T}         \text{.}
\]
We notice that the order of the elements~\( (x_{1}, \dotsc, x_{n}) \) does not matter.
In the case of a triple~\( (x, y, z) \), we get
\[
    \delta^{2} (x, y, z)
    =
    \diff{\delta} (y z ) \diff{\delta \differentialdeltaspace} x     +
    \diff{\delta} (x z ) \diff{\delta \differentialdeltaspace} y     +
    \diff{\delta} (x y ) \diff{\delta \differentialdeltaspace} z     +
    \diff{\delta \differentialdeltaspace} x \diff{\delta \differentialdeltaspace} y \diff{\delta \differentialdeltaspace} z     \text{.}
\]

We claim:

\begin{theorem}\label{res:generatorss_and_relations_stable_case}
    The \( \mathit{PS} \)-algebra~\( S \) is the free commutative algebra~\( F \)
    modulo the following relations:
    Firstly, for any odd~\( k \), the
    element~\( p_{k}^{2} = p_{2k } \) comes with the relation
    \begin{equation}\label{eq:p^2_relation}
		p_{k}^{2}
		=
		\sum_{i =0}^{(k -1)/2} \bigl[\diff{\delta} (w_{2i } w_{2(k -i) } ) +\diff{\delta \differentialdeltaspace} w_{2i } \diff{\delta \differentialdeltaspace} w_{2(k -i) } \bigr]     \text{.}
    \end{equation}
    Secondly, for any triple~\( (x, y, z) \) of square-free monomials in the variables~\( w_{2i } \),
    we have the relation
    \[
        \delta^{2} (x, y, z) = 0
        \text{.}
    \]
    This last relation, when written out, is equivalent to
    \begin{equation}\label{eq:delta^2(xyz)_relation}
        \diff{\delta} (y z ) \diff{\delta \differentialdeltaspace} x         +
	    \diff{\delta} (x z ) \diff{\delta \differentialdeltaspace} y 	    +
	    \diff{\delta} (x y ) \diff{\delta \differentialdeltaspace} z         =
        \diff{\delta \differentialdeltaspace} x \diff{\delta \differentialdeltaspace} y \diff{\delta \differentialdeltaspace} z         \text{.}
    \end{equation}
\end{theorem}

\begin{example}
	The first non-trivial relation of the form~\(
		\delta^{2} (x, y, z) = 0
	\)
	is
	\[
		\delta^{2} (w_{2} w_{4} , w_{2}, w_{4})
		=
		0
		\text{.}
	\]
	Using the standard forms
	\begin{align*}
		w_{2}^{2} w_{4} 			&=
			w_{2} w_{4} \diff{d\differentialdspace} w_{2} +w_{4} \diff{d\differentialdspace} w_{2} ^{2} +w_{4} \diff{d\differentialdspace} w_{4} 	\\
	\shortintertext{and}
		w_{2} w_{4}^{2} 			&=
			w_{2} \diff{d\differentialdspace} w_{2} ^{4} +w_{2} w_{6} \diff{d\differentialdspace} w_{2} +w_{2} w_{4} \diff{d\differentialdspace} w_{4} 	\\
			&\qquad
			+w_{2} \diff{d\differentialdspace} w_{2} \diff{d\differentialdspace} w_{6} +\diff{d\differentialdspace} w_{2} ^{2} \diff{d\differentialdspace} w_{6} +\diff{d\differentialdspace} w_{4} \diff{d\differentialdspace} w_{6} +w_{2} \diff{d\differentialdspace} w_{8} 	\text{,}
	\end{align*}
	we can write out this relation as
	\begin{align*}
		&\diff{\delta} (w_{2} w_{4} )^{2} +\diff{\delta} (w_{2} w_{4} ) \diff{\delta \differentialdeltaspace} w_{2} \diff{\delta \differentialdeltaspace} w_{4} +\diff{\delta} (w_{2} w_{6} ) \diff{\delta \differentialdeltaspace} w_{2} 	\\
		&\qquad
		+\diff{\delta \differentialdeltaspace} w_{2} ^{6} +\diff{\delta \differentialdeltaspace} w_{2} ^{3} \diff{\delta \differentialdeltaspace} w_{6} +\diff{\delta \differentialdeltaspace} w_{2} ^{2} \diff{\delta \differentialdeltaspace} w_{4} ^{2} +\diff{\delta \differentialdeltaspace} w_{2} ^{2} \diff{\delta \differentialdeltaspace} w_{8} +\diff{\delta \differentialdeltaspace} w_{4} ^{3} 		=
		0
		\text{.}
	\end{align*}
	This relation will become crucial later in \zcref{ex:n=4}
	(see also \zcref{rem:need_non-disjoint_xyz}).
\end{example}

\begin{remark}\label{rem:minimality_stable_case}
	The above presentation is not quite minimal.
	Indeed, for every odd~\( k \),
	the relation~\eqref{eq:p^2_relation} makes the generator~\( \diff{\delta \differentialdeltaspace} w_{2k } \), and hence the relation itself, redundant.
	However, after removing those, the remaining generators constitute a minimal set.
	To see this, note that in the relations~\(
		\delta^{2} (x, y, z) = 0
	\),
	all terms are products of two or more factors.
	Therefore, letting~\(
		S_{+} \subset S
	\)
	denote the augmentation ideal,
	these relations become zero in the vector space~\(
		S_{+}/S_{+}^{2}
	\).
	On the other hand, for odd~\( k \), the left-hand side of the relation~\eqref{eq:p^2_relation} also becomes zero,
	making the right-hand side a linear relation involving the generator~\(
		\diff{\delta \differentialdeltaspace} w_{2k} 	\).
	Therefore, the remaining generators form a basis for~\(
		S_{+}/S_{+}^{2}
	\),
	proving that they are indeed a minimal set of generators.
	
	With regard to the relations, the proof of \zcref{res:d(x dy) = dx dy_special_case} below
	shows that we can in fact limit the relations~\eqref{eq:p^2_relation}
	to the cases where \( z \)~has the form~\( z = w_{2i } \).
	We do not currently know whether further reductions are possible.
	
	Despite these generators and relations being redundant, we have chosen to keep them, as doing so makes both this presentation, as well as the derived presentation in \zcref{res:generators_and_relations_instable_case}, easier to work with.
\end{remark}

\begin{remark}\label{rem:need_non-disjoint_xyz}
	Note that in the relations~\(
		\delta^{2} (x, y, z) = 0
	\),
	where \( x \), \( y \), and~\( z \) are square-free
	monomials in the variables~\( w_{2i } \),
	the products
	\(
		x y 	\),
	\(
		x z 	\),
	and~\(
		y z 	\)
	need not be square-free.
	We stress that it is \emph{not} possible to do away with the relations where these products contain squares.
	Indeed, in \zcref{ex:n=4} below,
	the relation~\(
		\delta^{2} (w_{2} w_{4} , w_{2}, w_{4})
		=
		0
	\)~will be absolutely necessary.
\end{remark}

\begin{remark}\label{rem:Lie_algebra_stable_case}
	Once the theorem has been proved, we can drop the \( \diff{\delta} \)~notation
	and simply write~\( \diff{\delta \differentialdeltaspace} x \) as~\( \diff{d\differentialdspace} x \).
	When that happens, the relation \eqref{eq:delta^2(xyz)_relation} will become the Jacobi identity for the Lie bracket
	\[
		\lbrack x, y \rbrack
		=
		x \diff{d\differentialdspace} y 		+
		y \diff{d\differentialdspace} x 		=
		x \, \widebar{y} 		+
		\widebar{x} \, y 		=
		\diff{d} (\widebar{x} \, y )
		=
		\diff{d} (x y )
		+
		\diff{d\differentialdspace} x \diff{d\differentialdspace} y 	\]
	on~\( R \).
	Interestingly, the other relations,~\eqref{eq:p^2_relation},
	can also be written in terms of this Lie bracket
	since the summand on the right-hand side is~\(
		\lbrack w_{2i }, w_{2(k -i) } \rbrack
	\).
\end{remark}

The remainder of this section will be devoted to the proof of this theorem.
Until then, we denote by~\( \tilde{S} \)
the algebra defined in the theorem.

\begin{lemma}\label{res:delta^2_relation_holds_for_arbitrary_x}
    We have~\( \delta^{2} (x, y, z) = 0 \) in~\( \tilde{S} \)
    when \( y \) and~\( z \) are square-free monomials in the variables~\( w_{2i } \) while \( x \)~is an arbitrary element of~\( R \).
\end{lemma}

\begin{proof}
    Use~\zcref{res:existence_of_the_standard_form}
    to write~\( x \) in the standard form as a linear combination of elements of the form~\(
    	p a \diff{d\differentialdspace} w_{2i_{1} } \dotsm \diff{d\differentialdspace} w_{2i_{n} }     \)
    where \( p \)~is a square-free monomial in the variables~\( p_{2i -1} \)
    while \( a \)~is a square-free monomial in the variables~\( w_{2i } \).
    Then due to the definition of~\( \diff{\delta} \),
    both~\( p \) and the factors~\( \diff{d\differentialdspace} w_{2i } \)
    can be moved out of~\( \diff{\delta} \)
    in the equation~\eqref{eq:delta^2(xyz)_relation}.
    Thus the relation reduces to the relations~\( \delta^{2} (a, y, z) = 0 \).
\end{proof}

\begin{lemma}\label{res:d(x dy) = dx dy_special_case}
    For any~\( x \in R \) and \( y \)~a square-free monomial in the variables~\( w_{2i } \), we have~\(
        \diff{\delta}(x \diff{d\differentialdspace} y )
        =
        \diff{\delta \differentialdeltaspace} x \diff{\delta \differentialdeltaspace} y     \)
    in~\( \tilde{S} \).
\end{lemma}

\begin{proof}
    For this, we use induction on the
    number of factors in~\( y \),
    starting from the case where \( y \)~only has one factor, i.e.~\( y = w_{2i } \), in which case
    the statement follows from the definition
    of~\( \diff{\delta} \).
    Next, suppose that \( y \)~splits as a product~\( y = u v \)
    of monomials of smaller length.
    \zcref[S]{res:delta^2_relation_holds_for_arbitrary_x}
    shows that~\( \delta^{2} (x, u, v) = 0 \),
    and hence we get
    \begin{align*}
        \diff{\delta}(x\diff{d\differentialdspace} y )
            &=
            \diff{\delta} (x (u \diff{d\differentialdspace} v +v \diff{d\differentialdspace} u +\diff{d\differentialdspace} u \diff{d\differentialdspace} v ) )
    \\
            &= \diff{\delta} (x u \diff{d\differentialdspace} v )
                + \diff{\delta} (x v \diff{d\differentialdspace} u )
                + \diff{\delta} (x \diff{d\differentialdspace} u \diff{d\differentialdspace} v )
    \\
            &= \diff{\delta} (x u ) \diff{\delta \differentialdeltaspace} v                 + \diff{\delta} (x v ) \diff{\delta \differentialdeltaspace} u                 + \diff{\delta \differentialdeltaspace} x \diff{\delta \differentialdeltaspace} u \diff{\delta \differentialdeltaspace} v     \\
            &= \diff{\delta \differentialdeltaspace} x \diff{\delta} (u v )                 = \diff{\delta \differentialdeltaspace} x \diff{\delta \differentialdeltaspace} y     \text{,}
    \end{align*}
    as claimed.
\end{proof}

\begin{lemma}\label{res:delta(p_l^2 x) = p_l^2 delta(x)}
    For~\( x \in R \) and \( k \)~odd,
    we have~\(
    	\diff{\delta} (p_{k}^{2} x )
    	=
    	p_{k}^{2} \diff{\delta \differentialdeltaspace} x     \)
    in~\( \tilde{S} \).
\end{lemma}

\begin{proof}
	Notice that this does \emph{not} follow immediately from the definition of~\( \diff{\delta} \),
	which a priori only gives us~\(
		\diff{\delta} (p_{k} x )
		=
		p_{k} \diff{\delta \differentialdeltaspace} x 	\)
	when \( p_{k} \)~does not occur in the standard form of~\( x \).
    But it does follow from \zcref{res:d(x dy) = dx dy_special_case}
    using \zcref{res:p^2_formula}
    and~\eqref{eq:p^2_relation}.
\end{proof}

\begin{lemma}\label{res:delta(p_l x) = p_l delta(x)}
    For~\( x \in R \) and \( k \)~odd,
    we have~\(
    	\diff{\delta} (p_{k} x )
    	=
    	p_{k} \diff{\delta \differentialdeltaspace} x     \)
    in~\( \tilde{S} \).
\end{lemma}

\begin{proof}
    Because of the definition of~\( \diff{\delta} \), we may as well assume that \( x \)~is a square-free monomial in the variables \( p_{2i -1} \) and~\( w_{2i } \).
    If \( p_{k} \)~does not divide~\( x \),
    the claim follows immediately from the definition of~\( \diff{\delta} \).
    On the other hand, if \( p_{k} \)~does divide~\( x \),
    write~\( x = p_{k} y \)
    for a square-free monomial~\( y \)
    not divisible by~\( p_{k} \).
    Then the conclusion follows from \zcref{res:delta(p_l^2 x) = p_l^2 delta(x)}
    and the case we already considered.
\end{proof}

\begin{lemma}\label{res:d(x dy) = dx dy}
    For arbitrary~\( x, y \in R \),
    we have~\(
    	\diff{\delta} (x \diff{d\differentialdspace} y )
    	=
    	\diff{\delta \differentialdeltaspace} x \diff{\delta \differentialdeltaspace} y     \)
    in~\( \tilde{S} \).
\end{lemma}

\begin{proof}
    By bringing~\( y \) on the standard form
    and moving out factors of~\( \diff{d\differentialdspace} w_{2i } \) and~\( p_{2i -1} \) (the latter using \zcref{res:delta(p_l x) = p_l delta(x)}),
    this reduces to \zcref{res:d(x dy) = dx dy_special_case}.
\end{proof}

\begin{lemma}\label{res:delta(dy)=0}
    We have~\( \diff{\delta} (\diff{d\differentialdspace} y ) = 0 \) in~\( \tilde{S} \) for all~\( y \in R \).
\end{lemma}

\begin{proof}
    Previous lemma with~\( x = 1 \).
\end{proof}

\begin{corollary}
    We have~\( \delta^{2} (x_{1}, \dotsc, x_{n}) = 0 \)
    in~\( \tilde{S} \)
    for arbitrary~\(
    	x_{1}, \dotsc, x_{n}     	\in
    	R
    \).
\end{corollary}

\begin{proof}
    \zcref[S]{res:delta(dy)=0} tells
    us that~\( 
    	\diff{\delta} (\diff{d} (x_{1} \dotsm x_{n} ))
    	=
    	0
    \).
    Then the statement follows by applying
    the Thick Leibniz Rule to~\(
    	\diff{d} (x_{1} \dotsm x_{n} )
	\)
    and using~\zcref{res:d(x dy) = dx dy}.
\end{proof}

\begin{corollary}\label{res:delta(S)=0}
	We have~\( \diff{\delta \differentialdeltaspace} x = 0 \) in~\( \tilde{S} \) for all~\( x \in S \).
\end{corollary}

\begin{proof}
	Because of \zcref{res:S/I_exterior_algebra},
	we can write~\( x = p + \diff{d\differentialdspace} y \)
	for some square-free polynomial~\( p \in \mathit{PS} \)
	and some~\( y \in R \).
	Therefore,
	\(
		\diff{\delta \differentialdeltaspace} x 		=
		\diff{\delta \differentialdeltaspace} p + \diff{\delta} (\diff{d\differentialdspace} y )
		=
		0
		\in
		\tilde{S}
	\),
	with~\( \diff{\delta \differentialdeltaspace} p \) being zero by construction of~\( \diff{\delta} \),
	and~\( \diff{\delta} (\diff{d\differentialdspace} y ) \) being zero due to \zcref{res:delta(dy)=0}.
\end{proof}

\begin{lemma}
    There is a surjective map of \( \mathit{PS} \)-algebras~\( \tilde{S} \to S \) sending \( \diff{\delta \differentialdeltaspace} x \)
    to~\( \diff{d\differentialdspace} x \) for all~\( x \in R \).
\end{lemma}

\begin{proof}
    We initially define the map on the generators~\(
    	\diff{\delta \differentialdeltaspace} a     \)
    for square-free monomials~\( a \) in the variables~\( w_{2i} \),
    noting that the map is well-defined since all relations are mapped to zero
    due to \zcref{res:p^2_formula} and equation~\eqref{eq:d^2_thick_Leibniz_relation_inspiration}.
    Next we consider elements of the form~\(
    	x
    	=
    	p a \diff{d\differentialdspace} w_{2i_{1} } \dotsm \diff{d\differentialdspace} w_{2i_{n} }     \)
    where \( p \)~is a square-free monomial in the variables~\( p_{2i -1} \)
    while \( a \)~is a square-free monomial in the variables~\( w_{2i } \).
    We note from \zcref{res:d(x dy) = dx dy}
    that~\(
    	\diff{\delta \differentialdeltaspace} x     	=
       	p \diff{\delta \differentialdeltaspace} a \diff{\delta \differentialdeltaspace} w_{2i_{1} } \dotsm \diff{\delta \differentialdeltaspace} w_{2i_{n} }     \),
    which is mapped to~\(
    	\diff{d\differentialdspace} x     	=
       	p \diff{d\differentialdspace} a \diff{d\differentialdspace} w_{2i_{1} } \dotsm \diff{d\differentialdspace} w_{2i_{n} }     \) because of the case we already considered.
    For arbitrary~\( x \in R \), recall from \zcref{res:existence_of_the_standard_form}
    that \( x \)~is a linear combination
    of terms of the form~\(
    	p a \diff{d\differentialdspace} w_{2i_{1} } \dotsm \diff{d\differentialdspace} w_{2i_{n} }     \)
    as before, so by linearity,
    \( \diff{\delta \differentialdeltaspace} x \)~is sent to~\( \diff{d\differentialdspace} x \).
    This also shows that we can hit any~\( \diff{d\differentialdspace} x \in I \),
    proving that the map is surjective.
\end{proof}

Now we aim to prove that the kernel of~\( \tilde{S} \to S \) is zero, so denote this kernel by~\( K \).
We shall think of elements of~\( K \) as expressions~\(
	G(\diff{\delta \differentialdeltaspace} x_{1} , \dotsc, \diff{\delta \differentialdeltaspace} x_{n} )
	\in
	\tilde{S}
\)
for some polynomial~\(
	G
	\in
	\mathit{PS} \lbrack t_{1}, \dotsc, t_{n} \rbrack
\),
where \( x_{1}, \dotsc, x_{n} \)~are square-free monomials in the variables~\( w_{2i } \),
such that the corresponding expression~\(
	G(\diff{d\differentialdspace} x_{1} , \dotsc, \diff{d\differentialdspace} x_{n} )
	\in
	S
\)
is zero.

\begin{lemma}\label{res:kernel_K_no_constant_term}
    An element of the kernel~\( K \) can be written
    as~\( G(\diff{\delta \differentialdeltaspace} x_{1} , \dotsc, \diff{\delta \differentialdeltaspace} x_{n} ) \),
    where \( x_{1}, \dotsc, x_{n} \)~are square-free monomials in the variables~\( w_{2i } \)
    and \(
    	G
    	\in
    	\mathit{PS} \lbrack t_{1}, \dotsc, t_{n} \rbrack
    \)~has no constant term.
\end{lemma}

\begin{proof}
    Let~\( G(\diff{\delta \differentialdeltaspace} x_{1} , \dotsc, \diff{\delta \differentialdeltaspace} x_{n} ) \in K \),
    where~\( G \in \mathit{PS} \lbrack t_{1}, \dotsc, t_{n} \rbrack \)
    is a polynomial with constant term~\( g_{0} \in \mathit{PS} \).
    Since \(
    	g_{0}
    	=
    	G(\diff{d\differentialdspace} x_{1} , \dotsc, \diff{d\differentialdspace} x_{n} )
    	+
    	g_{0}
    	\in
    	S
    \)~is a sum of elements from~\( I \),
    we have~\( g_{0} \in I \).
    Also, as the map~\( \mathit{PS}/Q(\mathit{PS}) \to S/I \) is an isomorphism (c.f.~\zcref{res:S/I_exterior_algebra}), this shows that~\(
    	g_{0}
    	\in
    	Q(\mathit{PS})
    \).
    Therefore, we can apply the first relation~\eqref{eq:p^2_relation}
    to rewrite~\( g_{0} \)
    in terms of the generators~\( \diff{\delta \differentialdeltaspace} x \)
    for square-free monomials~\( x \) in the variables~\( w_{2i } \).
    In other words, the chosen polynomial expression~\( G(\diff{\delta \differentialdeltaspace} x_{1} , \dotsc, \diff{\delta \differentialdeltaspace} x_{n} ) \in \tilde{S} \)
    is equal to another polynomial expression in the generators of~\( \tilde{S} \), but with no constant term.
\end{proof}

\begin{proof}[Proof of \zcref{res:generatorss_and_relations_stable_case}]
	\zcref[S]{res:kernel_K_no_constant_term},
	combined with \zcref{res:delta(p_l x) = p_l delta(x),res:d(x dy) = dx dy}, shows that any element of~\( K \) can be written as~\(
		\diff{\delta \differentialdeltaspace} x 	\)
	for some~\( x \in R \).
	Furthermore, as this~\( \diff{\delta \differentialdeltaspace} x \)~lies in the kernel~\( K \)
	of~\( \tilde{S} \to S \), 
	we have~\(
		\diff{d\differentialdspace} x = 0 \in S
	\),
	and hence~\( x \in S \).
	By \zcref{res:delta(S)=0}, this means that~\(
		\diff{\delta \differentialdeltaspace} x 		=
		0
		\in
		\tilde{S}
	\).
	We conclude that~\( K = 0 \), which proves the theorem.
\end{proof}

\newpage

\chapter{The finite Grassmannian}\label{chap:finite_grassmannian}

As mentioned in the introduction, it is a classical fact that the \( \mathbb{F}_{2} \)-cohomology
of the finite real Grassmannian~\(
	\mathup{Gr} (n, m)
\)
is given by
\begin{equation}\label{eq:cohomology_finite_grassmannian}
	H^{\SemantexBullet} (\mathup{Gr} (n, m);\mathbb{F}_{2})
	=
	\mathbb{F}_{2} \lbrack w_{1}, \dotsc, w_{n} \rbrack /(\widebar{w}_{m +1}, \dotsc, \widebar{w}_{m +n})
	\text{,}
\end{equation}
where \( \widebar{w}_{k} \)~is interpreted using the same formulae
as in the ring~\( R \), but with the convention that~\( w_{j} = 0 \) for~\( j > n \).
In terms of Young diagrams, this can be represented by the ring of diagrams
contained in an \( n \times m \)~rectangle.
For brevity, we shall write~\(
	R_{n, m}
	=
	H^{\SemantexBullet} (\mathup{Gr} (n, m);\mathbb{F}_{2})
\).

\begin{proposition}\label{res:coho_of_finite_grassmannian_in_mixed_coords}
	The cohomology ring~\(
		R_{n, m}
	\)
	can be given a presentation in \enquote{mixed coordinates}
	(c.f.~\zcref{def:mixed_coordinates})
	as follows:
	As generators, we take~\(
		w_{i}
	\)
	for even numbers~\(
		i
	\)
	with~\(
		2
		\le
		i
		\le
		n
	\)
	and~\(
		p_{j}
	\)
	for odd numbers~\(
		j
	\)
	with~\(
		1
		\le
		j
		\le
		n
	\).
	The relations are~\( w_{k} = 0 \)
	for odd numbers~\( k \) with~\( n < k \le m + n \)
	and~\( \widebar{w}_{k} = 0 \)
	for all numbers~\( k \) with~\( m < k \le m + n \).
	Here, \( \widebar{w}_{k} \)~is interpreted recursively using the formula~\eqref{eq:omega_definition},
	and for odd~\( k \), \( w_{k} \)~is interpreted
	recursively using the formula~\eqref{eq:newton_identity}.
	Furthermore, for odd numbers~\( j \) with~\( n < j < m + n \),
	\( p_{j} \)~is interpreted using \eqref{eq:newton_identity}.
	In all three cases, we use the conventions that~\( w_{i} = 0 \) for even numbers~\( i > n \) and~\( p_{j} = 0 \) for odd numbers~\( j \ge m + n \).
\end{proposition}

\begin{proof}
	We see that all of these relations are satisfied in the ring~\(
		R_{n, m}
	\).
	Indeed, to see that~\( p_{j} = 0 \) for all odd~\( j \ge m + n \),
	we first use
	\eqref{eq:omega_definition} inductively to show
	that~\(
		\widebar{w}_{i} = 0
		\in
		R_{n, m}
	\)
	for all~\(
		i > m + n
	\).
	For~\( j > m + n \), any of the two formulae from \zcref{res:odd_and_even_formulae} then implies that~\( p_{j} = 0 \).
	If \( m + n \)~is odd, there is also the case~\( j = m + n \),
	which follows by choosing the formula from \zcref{res:odd_and_even_formulae}
	that does not include the term~\( w_{n} \widebar{w}_{m} \).
	
	To show that these relations are exhaustive,
	we denote by~\smash{\( \tilde{R}_{n, m} \)}
	the ring constructed in the proposition
	and notice that the only remaining relation from~\eqref{eq:cohomology_finite_grassmannian}
	that we still need to verify is~\smash{\(
		w_{k} = 0 \in \tilde{R}_{n, m}
	\)}
	for \emph{odd} numbers~\( k > m + n \).
	To this end,
	the assumption that~\( p_{j} = 0 \) for~\( j > m + n \)
	allows us to apply the \enquote{even} variant of \zcref{res:odd_and_even_formulae}
	inductively to see that~\smash{\(
		\widebar{w}_{i}
		=
		0
		\in
		\tilde{R}_{n, m}
	\)}
	for all odd~\( i > m + n \).
	Summarizing, we have
	\begin{align*}
		w_{i}
			& = 0
			\qquad
			\text{%
				for all~\( i \) with~\( n < i \le m + n \)
				and even~\( i > m + n \),%
			}
	\\
	\shortintertext{and}
		\widebar{w}_{i}
			& = 0
			\qquad
			\text{%
				for all~\( i \) with~\( m < i \le m + n \)
				and odd~\( i > m + n \).%
			}
	\end{align*}
	We now claim that
	\[
		w_{i} = 0
		\text{ for odd } i > m + n
		\qquad\text{and}\qquad
		\widebar{w}_{i} = 0
		\text{ for even } i > m + n
		\text{.}
	\]
	Assuming inductively that this claim holds for~\( i \)
	with~\( m + n < i < k \),
	\eqref{eq:omega_definition}~shows that
	\[
		0
		=
		\sum_{i =0}^{k} w_{i} \, \widebar{w}_{k -i} 		=
		w_{k} + \widebar{w}_{k}
		=
		\begin{cases*}
			w_{k}
				& if \( k \)~is odd,
		\\
			\widebar{w}_{k}
				& if \( k \)~is even.
		\end{cases*}
	\]
	This proves the claim for~\( k \) as well.
\end{proof}

\newpage

\begin{remark}
	If \( m + n \)~is odd, one of the relations in
	\zcref{res:coho_of_finite_grassmannian_in_mixed_coords} becomes redundant.
	Indeed, if \( n \)~is even and \( m \)~odd,
	the fact that~\( p_{m +n} = 0 \),
	together with the \enquote{odd} variant of
	\zcref{res:coho_of_finite_grassmannian_in_mixed_coords},
	implies that~\( w_{m +n} = 0 \).
	Similarly, if \( n \)~is odd and \( m \)~even,
	the same fact, together with the \enquote{even} variant,
	implies that~\( \widebar{w}_{m +n} = 0 \).
\end{remark}

From now on, we shall concentrate on the case~\( m = n \)
and write simply~\( R_{n} \) for the ring~\(
	R_{n, n}
	=
	H^{\SemantexBullet} (\mathup{Gr} (n,  n) ;\mathbb{F}_{2})
\),
represented by the Young diagrams contained in an \( n \times n \)~square.
We aim to formulate a version of the presentation in \zcref{res:generatorss_and_relations_stable_case} for
the subring~\(
	S_{n}
	\subset
	R_{n}
\)
of \( \omega \)-invariants.
Analogously to the stable case,
we denote by~\( I_{n} \subset S_{n} \)
the ideal consisting of elements of the
form~\(
	\diff{d\differentialdspace} x = x + \widebar{x}
\).

The first step towards such a presentation
is discarding the generators
from \zcref{res:generatorss_and_relations_stable_case}
which are clearly zero in~\( R_{n} \).
According to \zcref{res:coho_of_finite_grassmannian_in_mixed_coords},
these are the generators~\( p_{j} \) for odd~\( j > 2n \)
as well as~\( \diff{\delta \differentialdeltaspace} a \) for square-free monomials~\( a \) in the variables~\( w_{2i} \)
which contain~\( w_{j} \) for some even~\( j > n \).
Therefore, in what follows,
we shall use the convention that such generators are read as zero whenever they appear in the relations.
As in \zcref{chap:generators_and_relations_for_S},
we shall also allow ourselves to write~\( \diff{\delta \differentialdeltaspace} x \in S_{n} \)
for arbitrary~\( x \in R \), interpreted according to those same conventions.
Finally, if \( k \)~is an even number, express it as~\( k = 2^{i} j \) with~\( j \)~odd and write~\smash{\( p_{k} = p_{j}^{2^{i}} \)}.

\begin{theorem}\label{res:generators_and_relations_instable_case}
	The ring~\( S_{n} = R_{n}^{\omega} \) of \( \omega \)-invariants is the \( \mathbb{F}_{2} \)-algebra generated by~\(
		p_{j}
	\)
	for odd~\( j \) with~\( 1 \le j < 2n \)
	as well as symbols~\(
		\diff{\delta \differentialdeltaspace} a 	\)
	for square-free monomials~\( a \) in the variables~\(
		w_{i}
	\)
	with~\( i \) even and~\( 2 \le i \le n \),
	modulo the following relations:
	First and foremost, we have the relations
	\eqref{eq:p^2_relation} for odd~\( k \) with~\( 1 \le k < 2n \) as well as~\eqref{eq:delta^2(xyz)_relation}
	for square-free monomials~\( x, y, z \) in the variables~\( w_{i} \) with~\( i \) even and~\( 2 \le i \le n \).
	Finally, for odd~\( k \) with~\( n < k < 2n \),
	we have the relation
	\begin{equation}\label{eq:relation_delta(xw_k)=0}
		\diff{\delta} (x \, w_{k} )
		=
		0
	\end{equation}
	for all square-free monomials~\( x \) (including~\( x = 1 \))
	in the variables~\( w_{i} \)
	with~\( i \) even and~\( 2 \le i \le n \).
\end{theorem}

\begin{proof}
	We first observe that \( S_{n} \)~is a quotient
	of the ring~\( S \),
	namely by the ideal generated by~\(
		s_{\lambda}
	\)
	for partitions~\( \lambda = \lambda^{\!\vee} \)
	not contained in the \( n \times n \)~square
	and~\( \diff{\delta \differentialdeltaspace} s_{\mu} \)
	for partitions~\smash{\( \mu \neq \mu^{\vee} \)}
	not contained in the \( n \times n \)~square.
	Clearly, the relations discussed above are all satisfied
	in this ring.
	To prove that the relations are exhaustive, we must verify
	that they cause the elements
	\( s_{\lambda} \) and~\( \diff{\delta \differentialdeltaspace} s_{\mu} \) from before to be zero.
	To this end, denote by~\smash{\( \tilde{S}_{n} \)}
	the ring defined in the proposition.
	
	In the case of~\( s_{\lambda} \)
	with~\( \lambda = \lambda^{\!\vee} \),
	one proves, as in \zcref{res:S/I_exterior_algebra}, that
	\(
		S_{n}/I_{n}
	\)~is
	the exterior algebra~\(
		S_{n}/I_{n}
		=
		\bigwedge (p_{1}, p_{3}, \dotsc, p_{2n -1})
	\).
	Let~\( p \in \mathit{PS} \) be the element with all monomials square-free
	such that~\(
		s_{\lambda}
		= 
		p
		\in
		S/I
		=
		\bigwedge (p_{1}, p_{3}, p_{5}, \dotsc)
	\)
	(it is not hard to see that \( p \)~will be just one monomial, but we shall not need this).
	Since~\(
		s_{\lambda}
		=
		p
		=
		0
		\in
		S_{n}/I_{n}
	\),
	all monomials in~\( p \) must involve~\( p_{j} \)
	for some odd~\( j > 2n \),
	and hence~\smash{\( p = 0 \in \tilde{S}_{n} \)}.
	Furthermore, due to~\eqref{eq:diagram_formula_for_p_k},
	this means that \( p \)~is the sum of Schur functions associated to partitions not contained in the \( n \times n \)~square.
	Therefore, \( s_{\lambda} + p \in I \)~must be the sum of~\( \diff{\delta \differentialdeltaspace} s_{\mu} \) for some partitions~\( \mu \neq \mu^{\vee} \)
	not contained in the \( n \times n \)~square.
	In conclusion, we have reduced the problem to the second case,
	which we shall now consider.
	
	To prove that~\( \diff{\delta \differentialdeltaspace} s_{\mu} = 0 \)
	for~\( \mu \neq \mu^{\vee} \) not contained in the \( n \times n \)~square,
	we prove the more general claim
	\begin{equation}\label{eq:claim_delta(x)=0}
		\diff{\delta \differentialdeltaspace} x 		=
		0
		\qquad
		\text{for all~\( x \in R \) with~\( x = 0 \in R_{n} \).}
	\end{equation}
	By \zcref{res:coho_of_finite_grassmannian_in_mixed_coords},
	proving the claim~\eqref{eq:claim_delta(x)=0} boils down to proving
	for all~\( y \in R \)
	that~\(
		\diff{\delta} (y \, w_{k} )
		=
		0
	\)
	for odd~\( k \) with~\( n < k < 2n \)
	and even~\( k > n \),
	and that~\(
		\diff{\delta} (y \, \widebar{w}_{k} )
		=
		0
	\)
	for all~\( k \) with~\( n < k \le 2n \).
	But due to the calculation
	\begin{equation}\label{eq:rewriting_delta(yw_k)}
		\diff{\delta} (y \, \widebar{w}_{k} )
		=
		\diff{\delta} (y(w_{k} +\diff{d\differentialdspace} w_{k} ) )
		=
		\diff{\delta} (y \, w_{k} )
		+
		\diff{\delta \differentialdeltaspace} y \diff{\delta \differentialdeltaspace} w_{k} 		\text{,}
	\end{equation}
	this all boils further down to proving for all~\( y \in R \) that
	\begin{equation}\label{eq:claim_delta(xw_k)=0}
		\diff{\delta} (y \, w_{k} )
		=
		0
		\qquad
		\text{for odd~\( k \) with~\( n < k \le 2n \) and even~\( k \) with~\( k > n \).}
	\end{equation}
	As usual, by the properties of~\( \diff{\delta} \),
	it is enough to prove the claim~\eqref{eq:claim_delta(xw_k)=0}
	when \( y \)~is a square-free monomial in the variables~\( w_{i} \) for even~\( i \) (crucially, we can \emph{not} quite limit ourselves to~\( 2 \le i \le n \)).
	
	Suppose first that \( k \)~is odd and that~\( n < k < 2n \). If \( y \)~is a monomial in only the variables~\( w_{i} \) for even~\( i \) with~\( 2 \le i \le n \), then \eqref{eq:claim_delta(xw_k)=0}~is the relation~\eqref{eq:relation_delta(xw_k)=0}.
	On the other hand, if \( y \)~involves \( w_{i} \) for even~\( i \) with~\( i > n \), then we can use~\eqref{eq:newton_identity}
	to reduce the claim~\eqref{eq:claim_delta(xw_k)=0} to smaller values of~\( k \) until we eventually hit an even~\( k \ge 0 \).
	If this~\( k \) still satisfies~\( k > n \),
	then we have reduced to the case where \( k \)~was even.
	If~\( k \le n \), there are two cases to consider:
	If the monomial~\( y \) does not contain~\( w_{k} \),
	then
	\(
		\diff{\delta} (y \, w_{k} )
	\)~is one of the generators that we have modded out by, so it is zero.
	If \( y \)~does contain~\( w_{k} \),
	we can then use \zcref{res:formula_for_w^2}
	to eliminate the square
	from~\(
		\diff{\delta} (y \, w_{k} )
	\)
	and arrive at one of the generators that we have modded out by.
	
	Next, in the case where \( k = 2j > n \)~is even,
	there are again two cases to consider:
	If the square-free monomial~\( y \) does not contain~\( w_{2j} \), then \(
		\diff{\delta} (y \, w_{2j} )
	\)~is one of the generators that we have modded out by,
	hence it is zero.
	On the other hand, if \( y \)~does contain~\( w_{2j} \),
	we may
	write~\( y = z \, w_{2j} \)
	for a square-free monomial~\( z \).
	Recalling from \zcref{res:omega_is_normal}
	that~\(
		w_{2j} \, \widebar{w}_{2j} 		\in
		I
	\),
	we can use the formula~\smash{\(
		w_{2j}^{2}
		=
		w_{2j} \diff{d\differentialdspace} w_{2j} 		+
		w_{2j} \, \widebar{w}_{2j} 	\)}
	to obtain
	\[
		\diff{\delta} (y \, w_{2j} )
		=
		\diff{\delta} (z \, w_{2j}^{2} )
			=
			\diff{\delta \differentialdeltaspace} w_{2j} \diff{\delta} (z \, w_{2j} ) 			+
			w_{2j} \, \widebar{w}_{2j} \diff{\delta \differentialdeltaspace} z 			\text{,}
	\]
	where \(
		w_{2j} \, \widebar{w}_{2j} 	\)~is understood as an element of~\smash{\( \tilde{S}_{n} \)}.
	It is therefore enough to prove that~\smash{\(
		w_{2j} \, \widebar{w}_{2j} 		=
		0
		\in
		\tilde{S}_{n}
	\)}.
	As in the proof of \zcref{res:omega_is_normal}, we have
	\[
		w_{2j} \, \widebar{w}_{2j} 			=
			\sum_{i =0}^{2j -1} \diff{\delta} (w_{i} \, \widebar{w}_{4j -i} ) 			=
			\sum_{i =0}^{2j -1} \;\;\bigl[\diff{\delta} (w_{i} \, w_{4j -i} ) +\diff{\delta \differentialdeltaspace} w_{i} \diff{\delta \differentialdeltaspace} w_{4j -i} \bigr] 		\text{.}
	\]
	For even values of~\( i \), the summand on the right
	is zero as
	\(
		\diff{\delta} (w_{i} \, w_{4j -i} )
	\)
	and~\(
		\diff{\delta \differentialdeltaspace} w_{4j -i} 	\)
	are two of the generators that we have modded out by.
	For odd~\( i \), we can apply~\eqref{eq:newton_identity}
	to~\( w_{i} \) to obtain
	\[
		\diff{\delta} (w_{i} \, \widebar{w}_{4j -i} )
		=
		\sum_{t =0}^{i -1} p_{i -t} \diff{\delta} (w_{t} \, \widebar{w}_{4j -i} ) 		\text{.}
	\]
	Since~\( 4j - i > k > n \),
	we have~\(
		w_{t} \, \widebar{w}_{4j -i} 		=
		0
		\in
		R_{n}
	\).
	Furthermore, we notice that~\(
		\deg (w_{t} \, \widebar{w}_{4j -i} )
		<
		4j 	\).
	Therefore, assuming inductively that
	the original claim~\eqref{eq:claim_delta(x)=0}
	holds when the element~\( x \) is homogeneous and satisfies~\(
		\deg (x) < 4j 	\),
	we conclude that~\(
		\diff{\delta} (w_{i} \, \widebar{w}_{4j -i} )
		=
		0
	\),
	and hence also that~\(
		w_{2j} \, \widebar{w}_{2j} 		=
		0
	\).
	This finishes the argument.
\end{proof}

\newpage

\begin{remark}\label{rem:minimality_instable_case}
	If the odd number~\( k \) satisfies~\(
		1 \le k \le n/2
	\),
	then as in \zcref{rem:minimality_stable_case},
	the relation~\eqref{eq:p^2_relation} makes the generator~\(
		\diff{\delta \differentialdeltaspace} w_{2k} 	\),
	and the relation itself,
	redundant.
	If~\( n/2 < k < n \), then the generator~\(
		\diff{\delta \differentialdeltaspace} w_{2k} 	\)
	no longer exists, and the relation instead makes one of the other generators in the expression redundant, along with the relation. The natural canonical choice is~\(
		\diff{\delta} (w_{k -1} w_{k +1} )
	\).
	In principle, we could have also chosen to make~\(
		\diff{\delta} (w_{k -1} w_{k +1} )
	\)
	redundant in the previous case, but that would hardly make the presentation simpler.
	Finally, if~\(
		n \le k < 2n 	\), the right-hand side of~\eqref{eq:p^2_relation} is zero, and we get the (non-redundant) relation~\(
		p_{k}^{2} = 0
	\).
	
	With the mentioned generators removed, we claim that the remaining ones form a minimal generating set.
	Indeed, the argument is more or less the same as in \zcref{rem:minimality_stable_case}, except that we now also have the additional relations~\eqref{eq:relation_delta(xw_k)=0}.
	By construction of~\( \diff{\delta} \), each term of the expression~\(
		\diff{\delta} (x \, w_{k} )
	\)
	contains a generator of the form~\( \diff{\delta \differentialdeltaspace} a \) as a factor.
	However, as~\(
		\diff{\delta} (x \, w_{k} )
	\)
	has odd degree, each term must also include a factor of the form~\( p_{k} \) with~\( k \)~odd.
	In conclusion, these relations become trivial in the vector space~\(
		(S_{n})_{+}/(S_{n})_{+}^{2}
	\),
	and the mentioned set of generators forms a basis.
	
	As in \zcref{rem:minimality_stable_case}, we can also limit ourselves to the relations~\(
		\delta^{2} (x, y, z) = 0
	\)
	where~\( z = w_{i} \) with~\( i \) even and~\(
		2 \le i \le n
	\).
	We do not currently know if further reductions of the number of relations are possible.
\end{remark}

\begin{example}[Example: \( n = 2 \)]
	In this case, we end up with just two generators, \( p_{1} \) and~\( p_{3} \), with the relations~\( p_{1}^{3} = 0 \),
	\( p_{3}^{2} = 0 \), and
	\(
		p_{1}^{2} p_{3} 		=
		0
	\).
	Due to \zcref{rem:minimality_instable_case},
	we do not need the generator~\( \diff{\delta \differentialdeltaspace} w_{2} \).
\end{example}

\begin{example}[Example: \( n = 3 \)]
	We end up with three generators, \( p_{1} \), \( p_{3} \),
	and~\( p_{5} \),
	subject to the relations~\(
		p_{3}^{2}
		=
		0
	\),
	\(
		p_{5}^{2}
		=
		0
	\),
	\(
		p_{1}^{2} p_{3} 		=
		p_{1}^{5}
	\),
	and~\(
		p_{1}^{2} p_{5} 		=
		p_{1}^{7}
	\).
	The first and third relations together imply
	that~\smash{\(
		p_{1}^{8}
		=
		0
	\)}.
\end{example}

\begin{example}[Example: \( n = 4 \)]\label{ex:n=4}
	We get the generators~\(
		p_{1}
	\),
	\(
		p_{3}
	\),
	\(
		p_{5}
	\),
	\(
		p_{7}
	\),
	and
	\(
		\diff{\delta \differentialdeltaspace} w_{4} 	\).
	We do not need the generator~\(
		\diff{\delta} (w_{2} w_{4} )
	\)
	due to \zcref{rem:minimality_instable_case}.
	The remaining relations can be reduced to
	\begin{align*}
		p_{5}^{2}
		&= 0
	\\
		p_{7}^{2}
		&= 0
	\\
		p_{1}^{5} +p_{1}^{2} p_{3} +p_{1} \diff{\delta \differentialdeltaspace} w_{4} 		&= 0
	\\
		p_{1}^{7} +p_{1}^{2} p_{5} +p_{3} \diff{\delta \differentialdeltaspace} w_{4} 		&= 0
	\\
		p_{1}^{4} p_{3} +p_{1}^{2} p_{5} +p_{1} p_{3}^{2} 		&= 0
	\\
		p_{1}^{2} p_{7} +p_{3}^{3} +p_{1}^{2} p_{3} \diff{\delta \differentialdeltaspace} w_{4} 		&= 0
	\\
		p_{1}^{3} p_{3}^{2} +p_{1}^{2} p_{7} +p_{5} \diff{\delta \differentialdeltaspace} w_{4} 		&= 0
	\\
		p_{1}^{5} p_{3}^{2} +p_{3}^{2} p_{5} +p_{7} \diff{\delta \differentialdeltaspace} w_{4} +p_{3} \diff{\delta \differentialdeltaspace} w_{4} ^{2} 		&= 0
	\\
		p_{3}^{2} p_{5} +p_{1} p_{3}^{2} \diff{\delta \differentialdeltaspace} w_{4} 		&= 0
	\\
		p_{1}^{2} p_{3} p_{7} +p_{1}^{4} \diff{\delta \differentialdeltaspace} w_{4} ^{2} +\diff{\delta \differentialdeltaspace} w_{4} ^{3} 		&= 0
	\\
		p_{3}^{2} p_{7} +p_{3}^{3} \diff{\delta \differentialdeltaspace} w_{4} +p_{1} \diff{\delta \differentialdeltaspace} w_{4} ^{3} 		&= 0
	\text{.}
	\end{align*}
	The relation of degree~\( 12 \) comes from~\(
		\delta^{2} (w_{2} w_{4} , w_{2}, w_{4})
		=
		0
	\).
	We notice that this is the only relation that can make the generator~\( \diff{\delta \differentialdeltaspace} w_{4} \) nilpotent
	(c.f.~\zcref{rem:need_non-disjoint_xyz}).
	Through an elaborate set of reductions, the above relations
	imply that \smash{\( p_{1}^{8} = 0 \)}, \( p_{3}^{4} = 0 \), and~\( \diff{\delta \differentialdeltaspace} w_{4} ^{4} = 0 \).
\end{example}

\begin{remark}
	As in \zcref{rem:Lie_algebra_stable_case}, the relations from \zcref{res:generators_and_relations_instable_case} can also be written in the language of Lie algebras.
	Indeed, the additional relations~\eqref{eq:relation_delta(xw_k)=0} correspond to modding out by the Lie ideal consisting of~\( \lbrack x, \widebar{w}_{k} \rbrack \) for all~\( x \in R_{n} \).
	Alternatively, because of the calculation~\eqref{eq:rewriting_delta(yw_k)},
	one can also use~\( \lbrack x, w_{k} \rbrack \) for all~\( x \in R_{n} \).
\end{remark}

\chapter*{Statements and declarations}

The author declares that he has no conflict of interest.

\printbibliography

\end{document}